\RequirePackage{etoolbox}
\csdef{input@path}{{style/}{graphics/}}
\documentclass[aos,MSNbibl,nameyear,seceqn,dvips]{arximspdf}
\usepackage{graphicx}

\doi{10.1214/14-AOS1298}
\volume{43}
\issue{2}
\pubyear{2015}
\firstpage{847}
\lastpage{877}
\docsubty{FLA}

\makeatletter

\newcommand{\rrvert}{\vert}

\newcommand{\llvert}{\vert}
\newtheorem{theorem}{Theorem}[section]
\newtheorem{lemma}{Lemma}[section]

\newproclaim{remark}{Remark}[section]
\newproclaim{assumption}{Assumption}[section]
\newtheorem{corollary}{Corollary}[section]
\newproclaim{example}{Example}[section]
\makeatother

\begin{document}
\begin{frontmatter}

\title{Testing for pure-jump processes for~high-frequency~data}
\runtitle{Testing for pure-jump processes}

\begin{aug}
\author[A]{\fnms{Xin-Bing}~\snm{Kong}\corref{}\ead[label=e1]{kongxblqh@gmail.com}\thanksref{T1}},
\author[B]{\fnms{Zhi}~\snm{Liu}\ead[label=e2]{liuzhi@umac.mo}\thanksref{T2}}
\and
\author[C]{\fnms{Bing-Yi}~\snm{Jing}\ead[label=e3]{majing@ust.hk}\thanksref{T3}}
\runauthor{X.-B. Kong, Z. Liu and B.-Y. Jing}
\thankstext{T1}{Supported by NSF China 11201080 and the Humanity and Social
Science Youth Foundation of Chinese Ministry of Education (12YJC910003).}
\thankstext{T2}{Supported by FDCT of Macau (No. 078/2012/A3 and No. 078/2013/A3)
and NSFC No.~11401607.}
\thankstext{T3}{Supported by Hong Kong RGC Grants HKUST6019/12P and
Hong Kong RGC6022/13P.}

\affiliation{Soochow University, University of Macau and
Hong Kong University of Science and Technology}
\address[A]{X.-B. Kong\\
CASER and School of Mathematics\\
Soochow University\\
Shizi Road, Soochow\\
P. R. China\\
\printead{e1}}
\address[B]{Z. Liu\\
Department of Mathematics\\
University of Macau \\
Macau\\
\printead{e2}}
\address[C]{B.-Y. Jing\\
Department of Mathematics \\
Hong Kong University of Science and Technology\\
Clear Water Bay\\
Hong Kong\\
\printead{e3}}
\end{aug}

%
\received{\smonth{4} \syear{2014}}
%
\revised{\smonth{12} \syear{2014}}

%
\begin{abstract}
Pure-jump processes have been increasingly popular
in modeling high-frequency financial data, partially due to their
versatility and flexibility. In the meantime, several statistical
tests have been proposed in the literature to check the validity of
using pure-jump models. However, these tests suffer from several
drawbacks, such as requiring rather stringent conditions and having
slow rates of convergence. In this paper, we propose a different
test to check whether the underlying process of high-frequency data
can be modeled by a pure-jump process. The new test is based on the
realized characteristic function, and enjoys a much faster
convergence rate of order $O(n^{1/2})$ (where~$n$ is the sample
size) versus the usual $o(n^{1/4})$ available for existing tests;
it is applicable much more generally than previous tests; for example, it
is robust to jumps of infinite variation and flexible modeling of
the diffusion component. Simulation studies justify our findings and
the test is also applied to some real high-frequency financial data.
\end{abstract}

%
\pdfkeywords{62M05, 62G20, 60J75, 60G20, Ito semimartingale, pure-jump process, integrated volatility, realized characteristic function}

\end{frontmatter}

\section{Introduction}

It\^{o}'s semimartingales are widely used in modeling the log prices of
an asset since they fit many stylized features of asset returns, and
in option pricing due to absence of arbitrage in efficient market.
Mathematically, they consist of two parts: a continuous local
martingale term and a pure-jump process with both big and small
jumps. It\^{o}'s semimartingale with a continuous local martingale is in
common use in the literature, for example, the Black and Scholes (\citeyear{BlaSch73}) model
(geometric Brownan motion), the \citet{Mer76} model and \citet{Kou02}
model (geometric Brownan motion plus finitely many jumps).

On the other hand, in recent years pure-jump processes have also
been accepted as an alternative model for log price processes or
even the latent spot volatility process to the classic models
mentioned earlier; see, for example,
Todorov and Tauchen (\citeyear{TodTau10}, \citeyear{TodTau14}) and
references therein. The idea behind the pure-jump modeling is that
small jumps can eliminate the need for a continuous martingale.
Pure-jump models are also very flexible. They include the normal
inverse Gaussian [\citet{Ryd97};
Barndorff-Nielsen (\citeyear{Bar97}, \citeyear{Bar98})],
the variance gamma [Madan, Carr and Chang (\citeyear{MadCarCha98})], the CGMY model of
Carr et al. (\citeyear{Caretal03N2}), the time-changed L\'{e}vy mdoels of Carr
et al. (\citeyear{Caretal03N1}), the non-Gaussian Ornstein--Uhlenbeck-based models of
\citet{BarShe01} and the L\'{e}vy-driven
continuous-time moving average (CARMA) models of Brockwell (\citeyear{autokey10})
for the stochastic volatility. Pure-jump models have been
extensively used for general option pricing [\citet{HuaWu04};
\citet{BroDet04}; Levendorski{\u\i} (\citeyear{Lev04}); \citet{Sch06};
\citet{25}] and for foreign exchange option pricing [Huang and Hung (\citeyear{HuaHun05}); \citet{DaaMad05}; \citet{CarWu07}]. Other
applications of pure-jump models include reliability theory [\citet{Dro86}], insurance valuation [\citet{Bal05}] and fianancial
equilibrium analysis [Madan (\citeyear{Mad06})].

Statistically, this forces us to reconsider the necessity of
including the local martingale part driven by Brownian motion in
modeling high-frequency data. This begs the following question: ``Is
it sufficient to model high frequency data by pure-jump process
alone,'' or equivalently, ``is it necessary to add a Brownian force
underlying the high frequency data?'' The answer to this question
serves as a model selection purpose. For more motivation and
explanation, we refer to \citet{AtSJac10} and Jing, Kong and Liu (\citeyear{JinKonLiu12}).

For ease of presentation, let $X_t$ be a semimartingale defined on
some filtered probability space $(\Omega, \mathcal{F}, P)$,
\[
X_t=X_0+\int^t_0b_s\,ds+
\int^t_0\sigma_s\,dW_s+X^d_t,
\]
where $X_0$ is the initial value, $\int^t_0b_s\,ds$ is the drift term
with $b_s$ being the time-varying drift coefficient which is an
optional and c\`{a}dl\`{a}g process, $\int^t_0\sigma_s\,dW_s$ is a
continuous local martingale with $\sigma_s$ being an adapted process
and $W_s$ a standard Brownian motion and the last term is a
pure-jump component with the jump activity index $\beta$ defined by
%
\begin{equation}
\beta=\inf\biggl\{r; \sum_{0\leq s\leq T}|
\Delta_sX|^r\leq\infty\biggr\},
\end{equation}
where $\Delta_sX=X_s-X_{s-}$; see \citet{AtSJac09} and
Jing et al. (\citeyear{Jinetal12}). Then the above question can be formulated
as a hypothesis testing problem as
%
\begin{equation}
\label{hypothses} H_0\dvtx \int^T_0
\sigma_s^2\,ds>0\quad \mbox{v.s.} \quad H_1\dvtx \int
^T_0\sigma_s^2\,ds=0,
\end{equation}
where $T$ is the time span of the high-frequency data.

The testing problem (\ref{hypothses}) was studied by several
authors. For instance, \citet{ConMan07}, A\"{\i}t-Shalia and
Jacod (\citeyear{AtSJac10}) used threshold power variation to construct their test
statistics. However, there are two main drawbacks with the threshold
power variation method:
\begin{itemize}
\item First, their tests require that $X^d$ be of
finite variation, which
rules out many interesting models in finance since empirical
evidences in some real data analysis show that the jumps are of
infinite variation; see, for example, \citet{AtSJac09} and
\citet{ZhaWu09}.

\item Second, their tests are not very powerful, even when $\beta$
($0\leq\beta<2$)
is close to~$0$. This is rather counterintuitive since
probabilistically the smaller the value of $\beta$ is, the farther
$X^d$ is from a continuous semimartingale.
\end{itemize}

Interestingly, \citet{TodTau11} invented a test based on the
point estimator of the JAI cleverly constructed as the smallest power
for which the realized power variation (without thresholding) does not
explode. Surprisingly, a test based on this estimator for the presence
of Brownian motion has the property that it has more power for lower
level of activity. However, since it is from realized power variation,
once more, one has to assume that $X^d$ is of finite variation when
$C_T=\int^T_0\sigma_s^2\,ds$ does not vanish in order to have available
central limit theorem. It is also worth noticing that \citet{TodTau14} test for presence of Brownian motion by checking whether
``devolatilized'' truncated returns are i.i.d. normal assuming finite
activity jumps present in the underlying log price processes.

Testing the existence of a nonvanishing continuous local martingale
is challenging when the jumps are of infinite variation. Jing, Kong and Liu (\citeyear{JinKonLiu12}) used the number of small increments to propose a
test, which mitigates the above-mentioned difficulties, and can
handle jumps of infinite variation. However, it still has the
following deficiencies:
\begin{itemize}
\item
First, the local volatility model is too restrictive. For example, it
does not even cover the Heston model under $H_0$.


\item
Second, the spot volatility of the continuous component is assumed
to be positive almost everywhere in time $t$. So if $H_0$ is
rejected, it is quite possible that the continuous component
vanishes only in certain subintervals, but is still present in other
subintervals;
see the simulation in Section~\ref{s4} for more illustration.
\end{itemize}

In this paper, we develop a novel test to (\ref{hypothses}) to
overcome the difficulties encountered in previous approaches. The
convergence rate of our new test under $H_0$ is of order $n^{-1/2}$
when the jump component is of infinite variation, which is faster
than that of all existing tests. The idea of the test is based on
the realized characteristic function, which was introduced in
\citet{TodTau12} to investigate the distributional
property of volatilities at different time points; see also Todorov, Tauchen and Grynkiv (\citeyear{TodTauGry11}) and \citet{26}. With observable i.i.d.
increments of a class of L\'{e}vy process with either finite activity
or infinite activity jumps, Chen, Delaigle and Hall (\citeyear{CheDelHal10}) proposed a regression
method based on the empirical characteristic function to estimate the
parameters of the drift, scale, stable index and the distribution of
the jump size of a compound Poisson process, while in our paper, we
assume a flexible It\^{o} semimartingale with stochastic volatility and
stochastic coefficient of jump measures, and assume that the time lag
of successive observations shrinks to $0$ (high-frequency data) rather
than fixed, as implicitly assumed in Chen, Delaigle and Hall (\citeyear{CheDelHal10}).
However, we remark that direct application of the realized
characteristic function does not work in testing (\ref{hypothses}),
and some other novel statistical techniques are needed.

The paper is organized as follows. In Section~\ref{s2}, we give some
assumptions and introduce our test statistics. Main results are
presented in Section~\ref{s3}. Section~\ref{s4} gives some
simulation studies and real data analysis. The main proofs are
postponed to the \hyperref[app]{Appendix}, and the proofs of some lemmas are provided
in the supplementary material [\citet{KonLiuJin14}] to this paper.

Throughout the paper, we assume that the available data set is
$\{X_{t_i};  0 \le i\le n\}$ which is discretely sampled from $X$,
and is equally spaced in the fixed interval $[0, T]$, that is,
$t_i=i\Delta_n$ with $\Delta_n=T/n$ for $0 \le i\le n$. Denote the
$j$th one-step increment by
\[
\Delta^n_j X=X_{t_j}-X_{t_{j-1}},\qquad 1
\le j\le n.
\]

\section{Methodology} \label{s2}

The key idea behind our test statistic is that the characteristic
function of the increments of the It\^{o}'s semimartingale is dominated
by the continuous local martingale part.

For illustration, let us take the following simple example:
\[
X_t=\sigma W_t+\gamma Y_t,
\]
where $\sigma\ge0$ is a constant spot volatility, $\gamma$ is some
constant and $Y_t$ is a symmetric $\beta$-stable process. Then the
logarithm of the characteristic function is
%
\begin{equation}
\label{ch1} \log{\psi_n(u)}\equiv\log{E\bigl[e^{\sqrt{-1}u\Delta^n_iX/\sqrt{\Delta_n}}
\bigr]} =-\tfrac{1}{2}\sigma^2u^2-|
\gamma|^{\beta}u^{\beta}\Delta_n^{1-\beta/2}.
\end{equation}
As $\Delta_n \to0$, the last term in (\ref{ch1}) induced by the
jump part decreases at a rate of $\Delta_n^{1-\beta/2}$. Note that
when $\beta<1$ (i.e., $Y_t$ is of finite variation), in the context
of estimating $\sigma$ (or its functionals), the bias caused by the
jump part is of negligible size $o(\Delta_n^{1/2})$. This implies
that an estimator of $\sigma_t$ (or its functionals) for a general
semimartingale based on the characteristic function would very
likely be robust to jumps of finite variation, which is confirmed in
\citet{TodTau12} and \citet{26}. On the
other hand, the problem becomes more challenging when $\beta>1$
since the last term in (\ref{ch1}) is no longer a negligible bias
term. In testing (\ref{hypothses}), under $H_0$, the right-hand side
of (\ref{ch1}) is a nonvanishing constant while under $H_1$ it is
almost zero. This is a major feature we will explore later to
differentiate the null and the alternative hypotheses.

We shall now introduce our test statistic. To start with, we split
the data into $m_n$ nonoverlapping blocks with each block length
equal to $2v_n$ consisting of $2k_n$ intervals of length $\Delta_n$,
where $k_n$ is some integer depending on $n$. Motivated by~(\ref{ch1}), and in view of $X_{t+s}-X_{t} \approx
\sigma_{t}(W_{t+s}-W_t)+\gamma^{+}_{t-}(Y^{+}_{t+s}-Y^{+}_t)+\gamma
^{-}_{t-}(Y^{-}_{t+s}-Y^{-}_t)$ where $Y^{\pm}$ are two independent
``stable like'' L\'{e}vy processes and $\gamma^{\pm}$ are two c\`{a}dl\`
{a}g processes that will be specified later in Assumption~\ref{ass1}.
When $s$ is close to $0$, we can estimate $\sigma^2_{2jv_n}$ ($0 \le
j \le m_n-1$) locally by
%
\begin{equation}
\label{spot} c^0_j(u)=-\frac{1}{u^2}\log{
\biggl(L^0_j(u)\vee\frac{1}{\sqrt{k_n}}\biggr)},
\end{equation}
where
%
\begin{eqnarray}
L^0_j(u)&=&\frac{1}{k_n-1}\sum
^{k_n-1}_{l=1} \cos\bigl(u\bigl(\Delta^n_{2jk_n+2l+1}X-
\Delta^n_{2jk_n+2l}X\bigr)/\Delta_n^{1/2}
\bigr).
\end{eqnarray}
Summing over $c^0_j(u)$ for all $j\leq m_n$ and properly normalizing
it, one easily gets an estimator of the integrated volatility
process,
\[
C_t \equiv\int^t_0
\sigma_s^2\,ds.
\]
\citet{26} introduced a bias-corrected estimator of
$C_t$ as
%
\begin{equation}
\label{IV} \hat{C}_0(u_n)=2v_n\sum
^{[t/(2v_n)]-1}_{j=0} \biggl(c_j^0(u_n)-
\frac
{1}{u_n^2(k_n-1)}\bigl(\sinh\bigl(u_n^2c_j^0(u_n)
\bigr)\bigr)^2 \biggr),
\end{equation}
and further showed that
\begin{equation}
\label{decom} \hat{C}_0(u_n)=C_t+A_0(u_n)^n_t+O_p
\bigl(\Delta_n^{1/2}\bigr),
\end{equation}
where
\[
A_0(u)^n_t=2u^{\beta-2}
\Delta_n^{1-\beta/2} \int^t_0
a_s \,ds
\]
with
$a_s=\chi(\beta)(|\gamma_s^{+}|^{\beta}+|\gamma_s^{-}|^{\beta})$ and
$\chi(\beta)=\int^{\infty}_0 y^{-\beta} \sin y \,dy$. Then a natural
test statistic which can differentiate the null and alternative
hypotheses is
\[
T_n^{\prime}\equiv\frac{\hat{C}_0(2u_n)-\hat{C}_0(u_n)}{\hat{C}_0(u_n)} \longrightarrow_p
\cases{ %
0 ,& \quad $\mbox{on }\{C_T>0\},$
\vspace*{2pt}\cr
2^{\beta-2}-1<0, & \quad $\mbox{on }\{C_T=0\}.$}
\]
The problem with $T_n^{\prime}$ is that no central limit theorem is
available as $\beta>1$, so that one cannot find the rejection region
when jumps are of infinite variation. We will fix this problem with
some manipulations to $T_n^{\prime}$ below.

To do this, we replace $\hat{C}_0(u)$ by a similarly defined
quantity. 
Let the $c^1_j$'s and $\hat{C}_1(u)$ be similarly defined as the
$c^0_j$'s and $\hat{C}_0(u)$ with
$\Delta^n_{2jk_n+2l+1}X-\Delta^n_{2jk_n+2l}X$ replaced by
$\Delta^n_{2jk_n+2l}X-\Delta^n_{2jk_n+2l-1}X$, for $l=1,\ldots,
k_n-1$. A seemingly
better test statistic is then
%
\begin{equation}
\label{stat2} \qquad T_n^{*}=\frac{(\hat{C}_0(2u_n)-\hat{C}_1(u_n))-(\hat{C}_0(2u_n)-\hat
{C}_0(u_n))}{\hat{C}_1(u_n)} =
\frac{\hat{C}_0(u_n)-\hat{C}_1(u_n)}{\hat{C}_1(u_n)},
\end{equation}
which works under $H_0$ because the numerator is equal to
%
\begin{eqnarray}
\label{tstar} && \bigl[\bigl(\hat{C}_0(2u_n)-C_t-A_{0}(2u_n)^n_t
\bigr)-\bigl(\hat {C}_1(u_n)-C_t-A_0(u_n)^n_t
\bigr)\bigr]
\nonumber
\\
&&\quad -\bigl[\bigl(\hat{C}_0(2u_n)-C_t-A_0(2u_n)^n_t
\bigr)-\bigl(\hat {C}_0(u_n)-C_t-A_0(u_n)^n_t
\bigr)\bigr]
\\
&&\qquad =O_p\bigl( \Delta_n^{1/2}
\bigr)-o_p\bigl( \Delta_n^{1/2} \bigr).\nonumber
\end{eqnarray}
The second term in (\ref{tstar}) is $o_p(\Delta_n^{1/2})$ since
$\hat{C}_0(2u)$ and $\hat{C}_0(u)$ are calculated
in the same way, except for using different arguments, and are
asymptotically perfectly correlated as
$u=u_n\rightarrow0$; see also (a) in Theorem~1 of \citet{26}. However, the first term in (\ref{tstar}) is
$O_p(\Delta_n^{1/2})$ since $\hat{C}_1(u_n)$ uses the data points
one grid after those in $\hat{C}_0(2u_n)$, which decreases the
overlap of the data and hence has lower dependency between the terms
with argument $2u_n$ and $u_n$; see Theorem~\ref{th2} below.

Although ${T_n^*}/{\Delta_n^{1/2}}$ is tight under $H_0$, it can be
close to zero with a large probability under $H_1$ since the signal
in the numerator is swept away in the bias correction. This causes
difficulty in successfully detecting pure-jump processes under
$H_1$ and hence results in a low power. This difficulty can be remedied
by adding a bias of order $o(\Delta_n^{1/2})$ onto the numerator of
$T^*_n$.

Our final test statistic is
%
\begin{equation}
\label{ft} T_n=\frac{\hat{C}_0(u_n)-\hat{C}_1(u_n)-\gamma_n\Delta_n^{1/2}}{\hat{C}_1(u_n)},
\end{equation}
where $\gamma_n$ is some chosen constant satisfying
$\gamma_n\rightarrow0$ of which the explicit form will be given in
Section~\ref{tuning}. It can be shown that
%
\begin{equation}
\label{stat} T_n/\Delta_n^{1/2} \cases{ %
 = O_p(1), & \quad $\mbox{on } \{C_T>0\},$
\vspace*{2pt}\cr
\rightarrow^P -\infty&\quad $\mbox{on }\{C_T=0\}.$}
\end{equation}
This means that $T_n/\Delta_n^{1/2}$ can be used to differentiate
$H_0$ and $H_1$.


\section{Main results}\label{s3}

\subsection{Model assumptions}

We need the following assumptions.

\begin{assumption}\label{ass1}
\[
X^d_t=\int^t_0
\gamma^{+}_{s-}\,dY^{+}_s+\int
^t_0\gamma ^{-}_{s-}\,dY^{-}_s+
\int^t_0\int_R\delta(s,
z)p(ds, dz),
\]
where $Y^{+}$ and $Y^{-}$ are two independent L\'{e}vy processes
with positive jumps and L\'{e}vy triplet equal to $(0, 0,
F^{\pm})$, $\gamma^{\pm}$ are two c\`{a}dl\`{a}g adapted processes
and $p$ is a Poisson random measure on $R_{+}\times R$ with
intensity $q(dt, dx)=dt\otimes dx$. We assume further that, for some
$\beta> 1>r$, the L\'{e}vy measure satisfies
\[
\biggl|\overline{F}^{\pm}(x)-\frac{1}{x^{\beta}}\biggr|=\biggl|F^{\pm}\bigl((x,
\infty)\bigr)-\frac{1}{x^{\beta}}\biggr|\leq g(x),\qquad x\in(0, 1],
\]
with $g(x)$ a decreasing function s.t.
$\int^1_0x^{r-1}g(x)\,dx<\infty$, and $|\delta(t, x)|^r\wedge1\leq
J(x)$ with $J(x)$ Lebesgue integrable on $R$.
\end{assumption}


\begin{assumption}\label{ass2}
$\sigma_t$ is an It\^{o} semimartingale of the form
\begin{eqnarray*}
\sigma_t&=&\sigma_0 +\int^t_0b^{\sigma}_s\,ds+
\int^t_0H^{\sigma
}_s\,dW_s+
\int^t_0H^{\prime\sigma}_s\,dW^{\prime}_s
\\
& &{} +\int_0^t\int_{\{|\delta^{\sigma}(s, x)|\leq
1\}}
\delta^{\sigma}(s, x) (p-q) (ds, dx)\\
&&{}+\int^t_0
\int_{\{|\delta^{\sigma}(s, x)|> 1\}}\delta^{\sigma}(s, x)p(ds, dx),
\end{eqnarray*}
where all the integrands are optional processes satisfying the
integrable condition in It\^{o}'s sense, and $q$ is the compensator of
$p$. Assume that $W$ and $W^{\prime}$ are two independent Brownian
motions that are further independent of $(p, Y^+, Y^-)$.
\end{assumption}

\begin{assumption}\label{ass3}
We have a sequence $\tau_n$ of stopping times increasing to infinity, a
sequence $a_n$ of numbers and a nonnegative Lebsgue-integrable function
$J$ on R, such that the processes $b$, $H^{\sigma}$, $\gamma^{\pm}$ are
c\`{a}dl\`{a}g adapted, the coefficients $\delta$, $\delta^{\sigma}$
are predictable, the processes $b^{\sigma}$, $H^{\prime\sigma}$ are
progressively measurable and
\begin{eqnarray*}
&&\hspace*{-4pt}t<\tau_n \Rightarrow\bigl|\delta(t, z)\bigr|^r\wedge1 \leq
a_nJ(z), \bigl|\delta ^{\sigma}(t, z)\bigr|^2\wedge1 \leq
a_nJ(z),
\\
&&\hspace*{-4pt} t<\tau_n, V=b,b^{\sigma}, H^{\sigma},
H^{\prime\sigma}, \gamma^{\pm} \Rightarrow|V_t|\leq
a_n,
\\
&&\hspace*{-5pt} V=b, H^{\sigma}, \gamma^{\pm} \Rightarrow\bigl|E(V_{(t+s)\wedge\tau
_n}-V_{t\wedge\tau_n}|
\mathcal{F}_t)\bigr|+E\bigl(\bigl|V_{(t+s)\wedge\tau
_n}-V_{t\wedge\tau_n}|^2\bigr|
\mathcal{F}_t\bigr)\leq a_ns\hspace*{-1pt}.
\end{eqnarray*}
\end{assumption}

Assumption~\ref{ass1} is the same as the Assumption (A) given in \citet{26}. It essentially states that $X^d$ can be decomposed
into two
components: active and less active jumps. Here, the first two
components are the stable-like jumps assumed to have the jump
activity index $\beta>1$. (This can be extended to cover the case
for $r<\beta\leq1$ with extra efforts and possibly more stringent
conditions. However, if we have a priori that
$\beta<1$, more straightforward tests will be possible.) Another
reason we restrict attention to $\beta>1$ is
because this is more interesting and challenging statistically. The
last term consists of jumps with finite variation (but possibly of
infinite activity) which is expected to disappear in a limiting sense
as inspired by the finding following (\ref{ch1}). In
\citet{AtSJac10}, it is assumed that $\beta<1$
since otherwise no asymptotic distribution theory could be used
under $H_0$ to calculate the rejection region.

Assumption~\ref{ass2} is a standard assumption in the literature
which allows for the ``leverage'' effect due to the common driving
forces in $X$ and $\sigma$. In Assumption~\ref{ass2}, the jumps of
$\sigma_t$ are assumed, without restriction, to be driven by the same
Poisson measure as $X$.

Assumption~\ref{ass3} is the same as the Assumption (B) in \citet{26} and a rather general assumption which is
satisfied by the multifactor stochastic volatility models that are
widely used in financial econometrics, for example, the popular affine jump
diffusion models in Duffie, Pan and Singleton (\citeyear{Du00}).
Assumptions \ref{ass2} and \ref{ass3} admit a rather general It\^{o}
semimartingale as the continuous part under $H_0$. As a comparison,
Jing, Kong and Liu (\citeyear{JinKonLiu12}) require that the volatility be of form
$\sigma(X_t)$, a smooth function of $X_t$ bounded away from 0. Hence
our assumptions on the continuous component is far less restrictive
than that in Jing, Kong and Liu (\citeyear{JinKonLiu12}).

\subsection{Main theorems}

We first state a central limit theorem for the joint distribution of
$(\hat{C}_0(u_n), \hat{C}_1(u_n))$.

\begin{theorem}\label{th1}
Suppose $k_n$, $u_n$, $\gamma_n$ and $\Delta_n$ satisfy
%
\begin{eqnarray}
\label{cond} k_n\Delta_n^{1/2}&\rightarrow& 0,\qquad
k_n\Delta_n^{1/2-\varepsilon}\rightarrow\infty,\qquad
u_n\rightarrow0,\qquad \sup_n\frac{k_n\Delta_n^{1/2}}{u_n^4}<
\infty,
\nonumber
\\[-8pt]
\\[-8pt]
\nonumber
 \gamma_n&\rightarrow& 0,
\end{eqnarray}
for any $\varepsilon>0$. Let $c_s=\sigma_s^2$. Then on the set $\{C_t>0\}$
we have
%
\begin{eqnarray}
\label{clt} &&\frac{1}{\Delta_n^{1/2}}\pmatrix{
\hat{C}_0(u_n)-A_0(u_n)^n_t-C_t
\cr
\hat{C}_1(u_n)-A_0(u_n)^n_t-C_t}
\nonumber
\\[-8pt]
\\[-8pt]
\nonumber
&&\qquad
 \rightarrow^{\mathcal{L}_s} 2 \pmatrix{\displaystyle \int^t_0c_s
\,d\tilde{W}_s &
\vspace*{2pt}\cr
\displaystyle \int^t_0c_s d\biggl(
\frac{1}{2}\tilde{W}_s+\sqrt{3}/2\tilde{W}^{\perp}
\biggr) },
\end{eqnarray}
where $\tilde{W}$ and $\tilde{W}^{\perp}$ are two mutually
independent standard Brownian motions defined on an extension of the
original probability space and are further independent of
$\mathcal{F}$, and $\mathcal{L}_s$ stands for stable convergence.
\end{theorem}

In Theorem~1 of \citet{26}, a similar multivariate
central limit theorem related to the bias corrected estimator of $C_t$
in (\ref{IV}) with distinct arguments was obtained. While in (3.10) and
(3.11) of Theorem~1 of their paper, the vector of component estimators
with distinct multiples of $u_n$ are formed by using the same way of
aggregating the high-frequency data,\vspace*{1pt} Theorem~\ref{th1} in our paper
considers a bivariate central limit theorem for $ (\hat{C}_0(u_n),
\hat{C}_1(u_n) )$, with $\hat{C}_0(u_n)$ collecting the
high-frequency data one lag after $\hat{C}_1(u_n)$. By simple
application of Theorem~\ref{th1} and the continuous mapping theorem, we
soon have the following null distribution of
$T_n$.

\begin{theorem}\label{th2}
Under the conditions in Theorem~\ref{th1}, we have in restriction to $\{
C_t>0\}$,
\[
\Delta_n^{-1/2} T_n\rightarrow^{\mathcal{L}_s}
G_t,
\]
where $G_t$ is a centered Gaussian process with conditional variance
$\kappa_t=\frac{4\int^t_0c_t^2\,dt}{C_t^2}$.
\end{theorem}

It follows from Theorem~\ref{th2} that the convergence rate of $T_n$
is of order $\Delta_n^{1/2}$, in contrast to
$\Delta_n^{3/4-\varpi/2}$ in Jing, Kong and Liu (\citeyear{JinKonLiu12}), where
$\varpi>\beta-1/2$ is some constant (practically $\varpi$ is taken
as $3/2$ since $\beta$ is usually unknown) or
$v_n^{\beta^{\prime}/2}$ in \citet{AtSJac10}, where
$\beta^{\prime}<1$ and $v_n$ satisfies
\[
v_n/\Delta_n^{\rho-}\rightarrow0,\qquad
v_n/\Delta_n^{\rho
+}\rightarrow\infty,\qquad 0<\rho-<
\rho+<1/2.
\]
%

Theorem~\ref{th2} is not directly applicable in determining the
rejection region since the conditional variance is unknown. The
denominator of the conditional variance can be consistently estimated
by $(\hat{C}_1(u_n))^2$, thanks to (\ref{decom}). Inspired by the
construction of $\hat{C}_k(u)$ ($k=0, 1$), we use the following linear
combination of sample variances to estimate the integral in the
numerator of $\kappa_T$. Define
%
\begin{equation}
\label{Ihat} \hat{I}_n \equiv \tfrac{1}2 (
\hat{I}_{n,0} + \hat{I}_{n,1}), 
\end{equation}
where
%
\begin{equation}
\hat{I}_{n, k}=2v_n\sum^{[{t}/{(2v_n)}]-1}_{j=0}
\biggl(c_j^{k}(u_n)-\frac{(\sinh(u_n^2c^{k}_j(u_n)))}{u_n^2(k_n-1)}
\biggr)^2,\qquad  k=0, 1.
\end{equation}

Now we have the following studentized central limit theorem.

\begin{theorem}\label{th3}
Let
$\hat{\kappa}_T={4\hat{I}_n}/{(\hat{C}_1(u_n))^2}$.
Then we have under the conditions in Theorem~\ref{th1}, in restriction
to $\{C_T>0\}$,
%
\begin{equation}
\label{normalize} \mathcal{T}_n\equiv\frac{1}{\Delta_n^{1/2}}
\frac{T_n}{
\hat{\kappa}_T^{1/2} } \equiv \frac{\hat{C}_0(u_n)-\hat{C}_1(u_n)-\gamma_n\Delta_n^{1/2}}{2
\hat{I}_n^{1/2} \Delta_n^{1/2} }\rightarrow^{\mathcal{L}_s}
\mathcal{N}(0, 1),
\end{equation}
where $\mathcal{N}(0, 1)$ is a standard normal random variable
independent of $\mathcal{F}$.
\end{theorem}

%

From Theorem~\ref{th3}, we can reject $H_0$ if
$\mathcal{T}_n<-z_{\alpha}$ where $P(\mathcal{N}(0,
1)>z_{\alpha})=\alpha$ for $\alpha\in(0,1)$. Now we state a result
on the convergence rate of $\mathcal{T}_n$ under $H_1$.

\begin{theorem}\label{th4}
Suppose Assumptions \ref{ass1} and \ref{ass3} hold,
$k_n\Delta_n^{1/2}\rightarrow0$,  $k_n\Delta_n^{1/2-\varepsilon}
\rightarrow\infty$ for any $\varepsilon>0$, $\sup_n
{k_n\Delta_n^{1/2}}/{u_n^4}<\infty$ and $u_n$ is bounded. Then
on the set $\{C_t=0, \int^t_0a_s\,ds\neq0\}$, we have
%
\begin{equation}
\label{h1} \hat{C}_0(u_n)-\hat{C}_1(u_n)=O_p
\bigl(u_n^{-2}\Delta_n^{1-{\beta
}/{(2(\beta+1-r))}}+
u_n^{\beta/2-2}\Delta_n^{1-\beta/4}\bigr)
\end{equation}
and
%
\begin{equation}
\label{h2} \hat{I}_n=4u_n^{2\beta-4}
\Delta_n^{2(1-\beta/2)}\int^t_0a_s^2\,ds+o_p
\bigl(u_n^{2\beta-4}\Delta_n^{2(1-\beta/2)}\bigr).
\end{equation}\vadjust{\goodbreak}
\end{theorem}

The following result concerning the size and power performance of
the test is a straightforward consequence of Theorems \ref{th3} and
\ref{th4}.

\begin{corollary} \label{cor1}
\textup{(1)} Under the conditions in Theorem~\ref{th1}, we have
$P(\mathcal{T}_n<-z_{\alpha}|\{C_T\neq 0\})\rightarrow\alpha$;

\begin{longlist}[(2)]
\item[(2)]
under the conditions in Theorem~\ref{th4}, if
\[
\gamma_n\bigl(u_n^2\Delta_n^{{\beta}/{(2(\beta+1-r))}-{1}/{2}}
+u_n^{2-\beta/2}\Delta_n^{{\beta}/{4}-{1}/{2}}\bigr)
\rightarrow \infty,
\]
we have
$P(\mathcal{T}_n<-z_{\alpha}|C_T=0, \int^T_0a_s\,ds\neq 0)\rightarrow1$.
\end{longlist}
%
\end{corollary}

\begin{remark}\label{rem3}
Corollary~\ref{cor1} shows that our new test achieves asymptotic
nominal level $\alpha$ and the asymptotic power 1. It follows from
the proof of Corollary~\ref{cor1} that $\mathcal{T}_n$ goes to
$-\infty$ with rate
$O_p(\gamma_n(\frac{u_n^2}{\Delta_n})^{(2-\beta)/2})$ under $H_1$ and
conditions in 2.
Thus the test becomes more powerful as $\beta$ gets closer to $0$,
which will be further confirmed by our simulation studies. This
overcomes the drawbacks of the test by \citet{AtSJac10}.
\end{remark}

\subsection{Choice of tuning parameters} \label{tuning}

We now study how to choose tuning parameters $k_n$, $u_n$ and
$\gamma_n$.
The major role of $k_n$ is to balance the bias and variance of
$\hat{C}_0(u_n)-C_t$ and $\hat{C}_1(u_n)-C_t$. The larger the $k_n$,
the smaller the bias and the larger the variance.
Hence we could choose $k_n=-c^{\prime}\Delta_n^{1/2}\log{\Delta_n}$ for
some constant $c^{\prime}>0$.

Now we turn to $u_n$. The rationale for letting $u_n\rightarrow0$
under $H_0$ is to guarantee the 
convergence in probability in (\ref{cv1}).
As in \citet{26}, we choose $u_n$ so that
$u_n^2\int^T_0c_s\,ds\rightarrow0$ by setting $u_n={c}
(\log{(1/\Delta_n)})^{-1/30}\times\break  {\mathit{BV}_T}^{-1/2}$, where
$\mathit{BV}_T=({\pi}/{2})\sum^{n-1}_{i=1}|\Delta^n_iX\Vert \Delta^n_{i+1}X|$ is
the bipower variation, which is a consistent estimator of
$\int^T_0c_s\,ds$.
Another advantage of such choice of $u_n$ is that it
would be enlarged under $H_1$, which in turn increases the power as
is seen from Corollary~\ref{cor1} and Remark~\ref{rem3}. Choosing an
optimal $c$ is quite hard. In order not to incur much approximation
error in (\ref{cv1}), we suggest to choose small $c$ when $n$ is
moderate, say $c=0.18$. Simulation studies where the data is
generated from a fitted model (no guarantee of good fitting
accuracy) assuming $H_0$ given in \citet{26} show that
choosing $c$ around $0.18$ would work well.

Finally, we look at $\gamma_n$. On the one hand, $\gamma_n$ should
be close to $0$ under $H_0$ in order not to produce a big bias for
$\mathcal{T}_n$; on the other hand, $\gamma_n$ should converge to
$0$ with a rate of
$u_n^{-2}\Delta_n^{1/2-\beta/2(\beta+1-r)}+u_n^{\beta/2-2}\Delta
_n^{1/2-\beta/4}$
so that the test has good power. This is easily achieved by setting
$\gamma_n=c^*/\log{({u_n^{2}}/{\Delta_n})}$ when $u_n$ is determined
by the aforementioned method. To be conservative, one can choose
small $c^*$ when $n$ is moderate, say $c^*=0.2$.\vadjust{\goodbreak}

\section{Numerical experiments} \label{s4}

\subsection{Simulation studies}

In this section, we conduct simulation studies to check the
performance of the new test and make comparisons with the test given
in Jing, Kong and Liu (\citeyear{JinKonLiu12}). We first consider the performance on
control of type~I error probability.
As in \citet{26}, we generate simulation data for
$5000$ times from the following stochastic volatility model:
%
\begin{eqnarray}
\label{dg} X_t&=& X_0+\int^t_0
\sqrt{c_s}\,dW_s+0.5 Y_t, \qquad 0\leq t\leq T,
\\
\label{dg1} c_t&=&c_0+\int^t_00.03(1.0-c_s)\,ds+0.15
\int^t_0\sqrt{c_s}\,dW^{\prime}_s,
\end{eqnarray}
for $0\leq t\leq3T/4$ and $c_t\equiv0$ if $3T/4\leq t\leq T$. In
order to incorporate the leverage effect, we set $\operatorname{corr}(dW,
dW^{\prime})=-0.5$. The parameters in the volatility dynamic are
specified by fitting actual financial data in the same reference
paper. The volatility $c_t$ is a square root diffusion process which
is widely used in financial applications. We tuned $k_n$, $u_n$ and
$\gamma_n$ as in Section~\ref{tuning} with $c=0.18$ and $c^*=0.2$.
We consider $n=1170, 2340, 4680$ which corresponds to sample the
data per $20, 10, 5$ seconds, respectively.
In the simulation, we let $T$ be one day consisting of 6.5 trading hours.

\begin{table}
\caption{Empirical sizes and the empirical powers of the new test; the
nominal level is $5\%$; $(n=1170, k_n=50)$; $(n=2340, k_n=78)$;
$(n=4680, k_n=100)$}\label{table1}
\begin{tabular*}{\textwidth}{@{\extracolsep{\fill}}lcccccc@{}}
\hline
& \multicolumn{3}{c}{\textbf{Empirical sizes}} & \multicolumn{3}{c@{}}{\textbf{Empirical power}} \\[-6pt]
& \multicolumn{3}{c}{\hrulefill} & \multicolumn{3}{c@{}}{\hrulefill} \\
\multicolumn{1}{@{}l}{$\bolds{\beta}$} & \multicolumn{1}{c}{$\bolds{n=1170}$} & \multicolumn{1}{c}{$\bolds{n=2340}$} &
\multicolumn{1}{c}{$\bolds{n=4680}$} & \multicolumn{1}{c}{$\bolds{n=1170}$} & \multicolumn{1}{c}{$\bolds{n=2340}$} & \multicolumn{1}{c@{}}{$\bolds{n=4680}$}\\
\hline
$1.0 $ & $0.0610 $ & $0.0586 $ & $0.0574 $ & $0.9988 $ & $0.9998 $ &
$1.0000 $\\
$1.1 $ & $0.0616 $ & $0.0624 $ & $0.0610 $ & $0.9984 $ & $0.9990 $ &
$1.0000 $\\
$1.2 $ & $0.0640 $ & $0.0635 $ & $0.0634 $ & $0.9936 $ & $0.9986 $ &
$0.9996 $\\
$1.3 $ & $0.0604 $ & $0.0601 $ & $0.0608 $ & $0.9596 $ & $0.9948 $ &
$0.9986 $\\
$1.4 $ & $0.0522 $ & $0.0616 $ & $0.0616 $ & $0.6508 $ & $0.8414 $ &
$0.9650 $\\
$1.5 $ & $0.0566 $ & $0.0624 $ & $0.0610 $ & $0.2902 $ & $0.3810 $ &
$0.5290 $\\
$1.6 $ & $0.0612 $ & $0.0514 $ & $0.0524 $ & $0.1328 $ & $0.1698 $ &
$0.2138 $\\
$1.7 $ & $0.0594 $ & $0.0624 $ & $0.0554 $ & $0.0942 $ & $0.1068 $ &
$0.1208 $\\
$1.8 $ & $0.0578 $ & $0.0550 $ & $0.0594 $ & $0.0776 $ & $0.0804 $ &
$0.0804 $\\
$1.9 $ & $0.0572 $ & $0.0568 $ & $0.0558 $ & $0.0748 $ & $0.0790 $ &
$0.0728 $\\
\hline
\end{tabular*}
\end{table}

\begin{figure}

\includegraphics{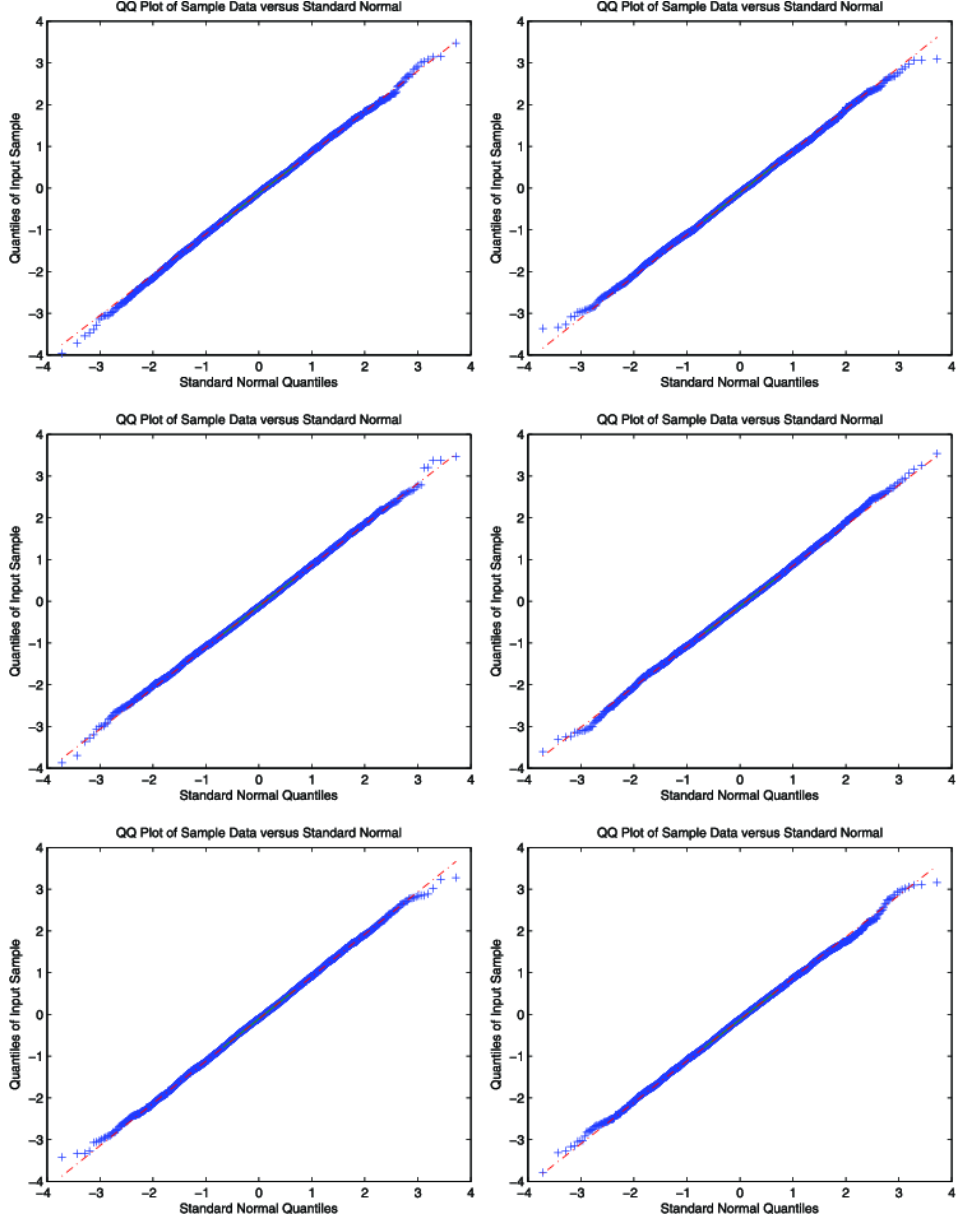}

\caption{QQ-plot of the test statistics under $H_0$ for
$\beta=1.2$ (left panels), 1.5 (right panels); from top to bottom,
$c=0.15, 0.18, 0.2$; $n=2340$.}\label{fig1}
\end{figure}

Table~\ref{table1} displays the empirical sizes of the new test.
Clearly, they are slightly higher than the nominal level but
acceptable across the board due to the small bias added
artificially.
Figure~\ref{fig1} gives the QQ-plot of the test statistics for
$n=2340$ and $\beta=1.2, 1.5$, showing that the normal
approximation works well.

For comparison, we choose $\beta=1.2$ and $n=2340$ and carry out
the test given in Jing, Kong and Liu (\citeyear{JinKonLiu12}), referred to as JKL's
test below. No comparisons will be made with the test given in
\citet{AtSJac10} (AJ's test), since it was
outperformed by the JKL's test in extensive simulation studies given
in Jing, Kong and Liu (\citeyear{JinKonLiu12}). Table~\ref{table2} lists the empirical
sizes of JKL's test where $\delta^*$ is a tuning parameter determining
how many small increments are used to compute the test statistics.
Clearly, the JKL's test is too liberal since the type I error
probabilities are out of control, showing that the JKL's test fails
when the continuous process vanishes in some subintervals. The
reason for the failure is that the JKL's test statistic has a
nonnegligible bias, even for large enough $n$.

\begin{table}
\caption{Empirical sizes of JKL's test; $\beta=1.2$, $n=2340$; the nominal
level is $5\%$; Empirical sizes$^*$ stand for the empirical sizes when
$c_s$ follows the same square root process for $3/4T\leq t\leq T$}
\label{table2}
\begin{tabular*}{\textwidth}{@{\extracolsep{\fill}}lcccccc@{}}
\hline
$\bolds{\delta^*}$ & \textbf{0.50} & \textbf{0.75} & \textbf{1.00} & \textbf{1.25} & \textbf{1.50} & \textbf{1.75}
\\
\hline
Empirical sizes & 0.3032 & 0.4442 & 0.5916 & 0.7312 & 0.8532 & 0.9402
\\
Empirical sizes$^*$ & 0.0298 & 0.0400 & 0.0358 & 0.0406 & 0.0402 &
0.0436 \\
\hline
\end{tabular*}
\end{table}

\begin{figure}[b]

\includegraphics{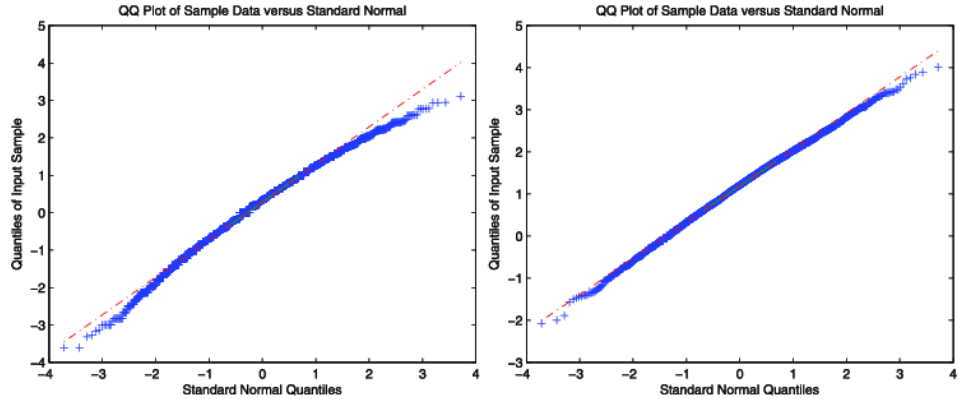}

\caption{QQ-plot of the test statistics given in \citet{Jinetal12} for $\delta^*=0.05$ (left panel),
0.5 (right panel) under $H_0$ when $\beta=1.2$; $n=2340$.}\label{fig2}
\end{figure}

It seems that choosing $\delta^*$ small would have satisfactory control
of type I error. However, when $\delta^*$ is small, the normal
approximation is actually no longer reliable. For $\delta^*=0.05$,
there are roughly $5$ small increments (effective data) used in
calculation of the test statistics, which affects the accuracy of the
normal approximation. Figure~\ref{fig2} gives the QQ-plot for the test
statistics given in Jing, Kong and Liu (\citeyear{JinKonLiu12})
for $\delta^*=0.05$ (left panel), 0.5 (right  panel) when
$\beta=1.2$ and $n=2340$. From
the left panel, we see a clear concavity pattern, which implies that
the distribution of the test statistic is left-skewed, yet the
empirical size is 0.07. Apparent improvement in skewness could be seen
in the right panel for $\delta^*=0.5$ since more effective data
(roughly 40) were added in calculation of the test statistics. However,
we see a clear bias in the QQ-plot.

Next we investigate the power of the new test. We generate the data
for $5000$ times from the above model, except that $c_s\equiv0$. The
empirical powers for various $\beta$ values are given in Table~\ref{table1}. We make the following observations:
\begin{itemize}
\item
the power of the new test decreases as $\beta$ increases since, as
$\beta$ increases to 2, the pure-jump process fluctuates more like a
Brownian motion;

\item
as the sample size increases, the empirical power increases overall,
as can be expected.
\end{itemize}

We also did a sensitivity study to $k_n$ when it is chosen in the
proposed range. In the sensitivity study we take $c^*=0.2$ and
$k_n=50, 78$, $c=0.15$ or $0.2$ when $n=2340$. The results on
both the size and power performance are reported in Table~\ref{table3}, where we can see that the empirical sizes and power do
not change much. We also conducted other sensitivity studies for $c
\approx0.18$ and $n=1170$ with $k_n$ in the corresponding range and
reached similar conclusions (hence not presented here).

\begin{table}
\caption{Empirical sizes and the empirical powers of the new test for
different pairs of $(c, k_n)$; the nominal level is $5\%$; n=2340}
\label{table3}
\begin{tabular*}{\textwidth}{@{\extracolsep{\fill}}lcccccc@{}}
\hline
&  \multicolumn{3}{c}{\textbf{Empirical sizes}} &
\multicolumn{3}{c@{}}{\textbf{Empirical power}} \\[-6pt]
&  \multicolumn{3}{c}{\hrulefill} &
\multicolumn{3}{c@{}}{\hrulefill} \\
\multicolumn{1}{@{}l}{$\bolds{\beta}$} & \multicolumn{1}{c}{$\bolds{(0.15, 50)}$} & \multicolumn{1}{c}{$\bolds{(0.15,  78)}$} &
\multicolumn{1}{c}{$\bolds{(0.2,  78)}$} & \multicolumn{1}{c}{$\bolds{(0.15, 50)}$} &
\multicolumn{1}{c}{$\bolds{(0.15, 78)}$} & \multicolumn{1}{c@{}}{$\bolds{(0.2,  78)}$}\\
\hline
$1.0 $ & $0.0630 $ & $0.0604 $ & $0.0634 $ & $0.9986 $ & $0.9994 $ &
$0.9992 $\\
$1.1 $ & $0.0604 $ & $0.0604 $ & $0.0608 $ & $0.9984 $ & $0.9986 $ &
$0.9992 $\\
$1.2 $ & $0.0634 $ & $0.0618 $ & $0.0624 $ & $0.9970 $ & $0.9982 $ &
$0.9982 $\\
$1.3 $ & $0.0558 $ & $0.0592 $ & $0.0638 $ & $0.9842 $ & $0.9896 $ &
$0.9968 $\\
$1.4 $ & $0.0580 $ & $0.0562 $ & $0.0618 $ & $0.7432 $ & $0.7708 $ &
$0.8786 $\\
$1.5 $ & $0.0584 $ & $0.0560 $ & $0.0614 $ & $0.3148 $ & $0.3242 $ &
$0.4146 $\\
$1.6 $ & $0.0576 $ & $0.0613 $ & $0.0608 $ & $0.1670 $ & $0.1498 $ &
$0.1814 $\\
$1.7 $ & $0.0558 $ & $0.0496 $ & $0.0568 $ & $0.0908 $ & $0.0906 $ &
$0.1102 $\\
$1.8 $ & $0.0558 $ & $0.0540 $ & $0.0582 $ & $0.0780 $ & $0.0778 $ &
$0.0788 $\\
$1.9 $ & $0.0542 $ & $0.0544 $ & $0.0566 $ & $0.0702 $ & $0.0702 $ &
$0.0744 $\\
\hline
\end{tabular*}
\end{table}

\subsection{Real data analysis}

In this section, we implement our test on some real data sets. We first
investigate
the stock price records of Microsoft (MFST) in two trading days,
December~1, and 12, 2000, which were also included in Jing, Kong and Liu
(\citeyear{JinKonLiu12}). All data sets are from the TAQ database. As in
Jing, Kong and Liu (\citeyear{JinKonLiu12}), to weaken the possible effect from microstructure
noise, we sample observations every $1/3$ minutes.
Finally, we use logarithms of the sampled prices to calculate the
test statistics.

\begin{figure}

\includegraphics{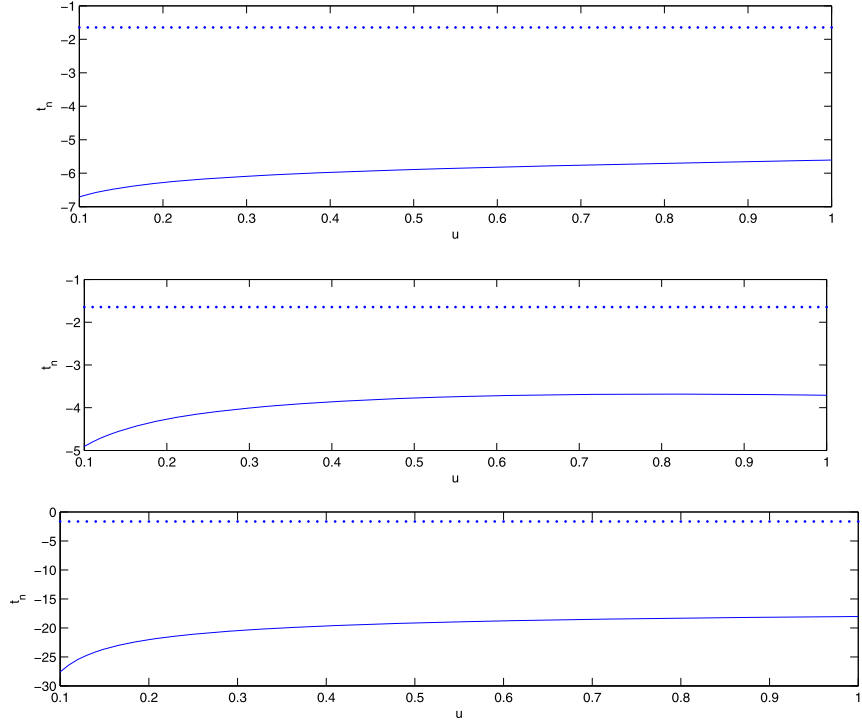}

\caption{Observed test statistics for the trading date, December 1 (middle
panel) and December~12 (upper panel),
in 2000, and the 5-mins S\&P 500 index (lower panel) data during January
4--29, 2010. The horizontal line has level $-1.645$.}\label{fig3}
\end{figure}

We set $T=1$ (day) consisting of 6.5 hours of trading time. As in
the simulation studies, we set $k_n=50$ and $\gamma_n=0.2/\log
{({u_n^{2}}/{\Delta_n})}$. To be on the safe side, let
$u_n$ take values in the grid points in $(0, 1]$ with step length
equal to $0.01$. Figure~\ref{fig3} plots the test statistics against
$u_n$ for two data sets. We see from the figure that for all
configurations of $u_n$, the test statistics are far lower than
$-1.645$, hence providing significant evidence against the existence of a
Brownian force. This confirms the empirical results in Jing, Kong and Liu (\citeyear{JinKonLiu12}) and in the meantime rules out the possibility that
Brownian force exists in some subintervals.

Next we implement our test the S\&P 500 index data which are sampled
every 5 minutes during January 4--29, 2010. The tuning parameters are used
as given above for those two stock data. The observed test statistics
are plotted against $u$ in the lower panel of Figure~\ref{fig3}. We
obtain the same conclusion that during the specified time period, the
underlying log price should be modeled by a pure-jump process.

\section{Conclusion and discussion}

In this paper, we have developed a new test based on the realized
characteristic function\vadjust{\goodbreak} to check whether the underlying process of a
high frequency data set can be modeled as a pure-jump process, and
shown its advantages over existing tests. Here are some future
problems worth pursing in future research work:
\begin{itemize}
\item
The effect of the microstructure noise, in the testing problem
(\ref{hypothses}) or even in estimating the functionals of the
volatility, is unclear and worthy of investigation in both
theory and practice. Here we could explore the two-time-scale
technique or multi-time-scale technique [A{\"{\i}}t-Sahalia, Mykland and Zhang (\citeyear{AtSMykZha05}),
\citet{Zha06}] or the pre-averaging approach [\citet{Jacetal09}].

\item
In the present paper, our inference is with the price process. It is
of interest to make inference on the volatility process which, as
recommended in \citet{TodTau14}, could be modeled by a
pure-jump process. The challenge of this problem is that the
volatility process is unobservable. Studies on this topic is still
undergoing.
\end{itemize}

\begin{appendix}
\section*{Appendix: Proofs of main theorems}\label{app}
This appendix contains the proofs of main theorems. The proofs of
Lemmas~\ref{lem3}--\ref{lem4} as well as some interesting supplemental lemmas are given
in Kong, Liu and Jing (\citeyear{KonLiuJin14}), a supplementary material [\citet{KonLiuJin14}] to this paper that is not for
purpose of publication. By the standard localization procedure, it is
enough to prove the main results under the following strengthened assumption.

\begin{assumption}\label{ass4}
$b$, $\sigma$, $\gamma^+$, $\gamma^-$, $b^{\sigma}$, $H^{\sigma}$ and
$H^{\prime\sigma}$ are bounded.
\end{assumption}

Before we prove the theorems, we introduce some notation and give
an outline of our proof. Let
$U_{t}(u)=\exp{(-u^2c_t-2\Delta_n^{1-\beta/2}u^{\beta}a_t)}$ where
$a_t=\chi(\beta)(|\gamma^+_t|^{\beta}+|\gamma_t^-|^{\beta})$ with
$\chi(\beta)=\int^{\infty}_0 y^{-\beta} \sin(y) \,dy$. For ease of
notation,\break $U_j(u)\equiv U_{2jv_n}(u)$ and sometimes we write
$E_{\mathcal{F}_t} V_s=E(V_s|\mathcal{F}_t)$ for a stochastic
process $V_t$. Let $\xi_{k,j}(u)=L_j^k(u)/U_j(u) - 1$, $k=0, 1$.
Let $\Omega(k, n, t)=\{\omega,\break \max_{k, j}|\xi_{k,j}(u,   \omega)|\leq
1/2\}$. By Lemma~7 of \citet{26},
%
\begin{equation}
\label{omega} P\bigl(\Omega^c(k, n, t)\bigr)\rightarrow0,
\end{equation}
irrespective of whether the continuous component exists or not.

\subsection{Proof of results under $H_0$}

Assuming the continuous local martingale exists, our proof depends
heavily on the following decomposition:
%
\begin{equation}\qquad\quad
\label{d1} c^k_j(u)=c_{2jv_n}+2u^{\beta-2}
\Delta_n^{1-\beta/2}a_{2jv_n}-\frac
{1}{u^2}
\xi_{k,j}(u)+\frac{1}{2u^2}\xi_{k,j}^2(u)+r_{k,j}(u),
\end{equation}
where $r_{k,j}(u)$ represents the remaining term which will be shown to
be negligible. By summing up the terms in (\ref{d1}) over $j$, one soon has
%
\begin{eqnarray}
\label{d2} \hat{C}_k(u)
&=&\sum^{[{t}/{(2v_n)}]-1}_{j=0}2v_nc_{2jv_n}+
\sum^{[
{t}/{(2v_n)}]-1}_{j=0}2v_n
\bigl(2u^{\beta-2}\Delta_n^{1-\beta/2}a_{2jv_n}\bigr)\nonumber\\
&&{}-
\sum^{[{t}/{(2v_n)}]-1}_{j=0}2v_n
\xi_{k,j}(u)/u^2
\\
& &{} +\sum^{[{t}/{(2v_n)}]-1}_{j=0}2v_n
\biggl(\frac{\xi_{k,j}^2(u)}{2u^2}-\frac
{1}{(k_n-1)u^2}\bigl(\sinh\bigl(u^2c^k_j(u)
\bigr)\bigr)^2\biggr)+R_k(u).\nonumber
\end{eqnarray}

We will first show that the first and second term converge to some
limits, and the fourth and last term in (\ref{d2}) are
$o_p(\Delta_n^{1/2})$, while the third term is $O_p(\Delta_n^{1/2})$
and converges to a conditionally centered Gaussian random variable
stably. This proves the univariate central limit theorem in Theorem~\ref{th1}. After that we proceed with the proof of the bivariate central
limit theorem by investigation into the covariation of those two
marginal sequences, which ends up with Theorem~\ref{th1}. Theorem~\ref{th2} is a consequence of Theorem~\ref{th1} and the continuous
mapping theorem. Theorem~\ref{th3} can be proved by showing that
$\hat{\kappa}_T$ is consistent to $\kappa_T$. In the sequel, $K$
will be a constant that has different values at different
appearances.

We now cite three lemmas from \citet{26}, whose proof can
be found in the same reference paper.
Lemma~\ref{lem1} is concerned with the first and second term in (\ref{d2}),
that is, the discretization error terms. Lemma~\ref{lem2} gives the stochastic
order of $\xi_{k,j}(u_n)$, $k=0, 1$, while Lemma~\ref{lem21} shows that the
fourth term and the remainder term in (\ref{d2}) are asymptotically
negligible.

\begin{lemma}[{[Lemma~8 in \citet{26}]}]\label{lem1}
 Under Assumptions~\ref
{ass1}--\ref{ass4} and assuming (\ref{cond}), we have
%
\begin{eqnarray}
\label{l11} \sum^{[t/(2v_n)]-1}_{j=0}2v_nc_{2jv_n}-
\int^t_0c_s\,ds&=&o_p
\bigl(u_n^2\Delta_n^{1/2}\bigr),
\\
\label{l12} \sum^{[t/(2v_n)]-1}_{j=0}2v_n
\bigl(2u_n^{\beta-2}\Delta_n^{1-\beta
/2}a_{2jv_n}
\bigr)-A_0(u_n)^n_t&=&o_p
\bigl(u_n^2\Delta_n^{1/2}\bigr).
\end{eqnarray}
\end{lemma}

\begin{lemma}[{[Lemma~14 in \citet{26}]}]\label{lem2}
 Under Assumptions~\ref
{ass1}--\ref{ass4} and assuming (\ref{cond}), we have, for $k=0, 1$,
%
\begin{eqnarray}
\label{l21} \bigl|E_{\mathcal{F}_{2jv_n}}\xi_{k,j}(u_n)\bigr|&\leq&
Ku_n^4\Delta_n^{1/2}
\phi_n,
\\
\label{l22} \biggl|E_{\mathcal{F}_{2jv_n}}\xi_{k,j}^2(u_n)-
\frac
{U_j(2u_n)+U_j(0)-2U_j^2(u_n)}{2(k_n-1)U_j^2(u_n)}\biggr|&\leq& Ku_n^4\Delta_n^{1/2}
\phi_n,
\end{eqnarray}
and for $q\geq2$,
%
\begin{equation}
\label{l23} E_{\mathcal{F}_{2jv_n}}\bigl|\xi_{k,j}(u_n)\bigr|^q
\leq K \bigl( u_n^{2q}/k_n^{q/2} +
u_n^4v_n \bigr),
\end{equation}
where $\phi_n$ is some sequence of numbers converging to 0.
\end{lemma}

\begin{lemma}[{[Lemma~9 in \citet{26}]}]\label{lem21}
 Under Assumptions~\ref
{ass1}--\ref{ass4} and assuming (\ref{cond}), we
have $R_k(u)=o_p(u_n^2\Delta_n^{1/2})$ $k=0, 1$ and
%
\begin{equation}\qquad
\label{l13} \sum^{[t/(2v_n)]-1}_{j=0}2v_n
\biggl(\frac{\xi_{k,j}^2(u_n)}{2u_n^2}-\frac
{1}{(k_n-1)u_n^2}\bigl(\sinh\bigl(u_n^2c^k_j(u_n)
\bigr)\bigr)^2\biggr)=o_p\bigl(u_n^2
\Delta_n^{1/2}\bigr).\vspace*{-6pt}
\end{equation}
\end{lemma}

The following lemma provides a formula for the limit of the conditional
real part of the characteristic function of a linear combination of
three successive increments. The proof can be found in the supplementary material [\citet{KonLiuJin14}]
to this paper.

\begin{lemma}\label{lem3}
Let $u_n^*=|a_{n, 0}|\vee|a_{n, 1}| \vee|a_{n, 2}|$, under
Assumptions \ref{ass1}--\ref{ass4}, and assume (\ref{cond}) with
$u^*_n$ replacing $u_n$, so we have
%
\begin{eqnarray}
\label{l31} && \Biggl|E_{\mathcal{F}_{(i-1)\Delta_n}} \cos\Biggl(\sum
_{l=0}^2a_{n,l}\frac
{\Delta^n_{i+l}X}{\Delta_n^{1/2}}\Biggr)
\nonumber
\\
&&\hspace*{10pt}\quad{} -\exp \Biggl(-\frac{1}{2}\sigma_{(i-1)\Delta_n}^2\sum
_{l=0}^2a_{n,
l}^2\nonumber\\
&&\hspace*{56pt}{}+
\Delta_n^{1-\beta/2}\chi(\beta)\sum^2_{l=0}
\bigl(\bigl|a_{n, l}\gamma ^+_{(i-1)\Delta_n}|^{\beta}+|a_{n, l}
\gamma^-_{(i-1)\Delta_n}\bigr|^{\beta
}\bigr) \Biggr)
\\
&&\hspace*{20pt}\quad{} \times\cos \Biggl(\Delta_n^{1-\beta/2}\chi^{\prime}(
\beta)\sum^2_{l=0}\bigl(\bigl
\{a_{n, l}\gamma^+_{(i-1)\Delta_n}\bigr\}^{\beta}+\bigl
\{a_{n, l}\gamma ^-_{(i-1)\Delta_n}\bigr\}^{\beta}\bigr) \Biggr) \Biggr|
\nonumber
\\
&&\qquad\leq Ku_n^{*4}\Delta_n^{1/2}
\phi_n,\nonumber
\end{eqnarray}
where $\{x\}^{\beta}=\operatorname{sign}{(x)}|x|^{\beta}$ and
$\chi^{\prime}(\beta)=\int^{\infty}_0\frac{1-\cos(y)}{y^{\beta}}\,dy$.
\end{lemma}

\begin{pf*}{Proof of Theorem~\ref{th1}}
By Lemmas \ref{lem1}, \ref{lem21} and (\ref{d2}), it suffices to prove that
\[
\frac{1}{\Delta_n^{1/2}} \Biggl(\sum^{[t/(2v_n)]-1}_{j=0}2v_n
\xi _{0,j}(u_{n})/u_{n}^2, \sum
^{[t/(2v_n)]-1}_{j=0}2v_n\xi
_{1,j}(u_{n})/u_{n}^2 \Biggr)
\]
converges to the right-hand side of (\ref{clt}) stably. By Lemma~\ref
{lem2}, we have
\[
\sum^{[t/(2v_n)]-1}_{j=0}2v_nE \bigl(
\xi_{k,j}(u_{n})/u_{n}^2|
\mathcal{F}_{2jv_n} \bigr)=o_p\bigl(u_n^2
\Delta_n^{1/2}\bigr),\qquad k=0, 1.
\]
Hence it is enough to prove the bivariate central limit theorem with
stable convergence for the following centered discrete bivariate
martingale with respect to $(\mathcal{F}_{2jv_n})_{j=0}^{[{t}/{(2v_n)}]-1}$:
%
\begin{eqnarray}
\label{bclt} && \frac{2v_n}{\Delta_n^{1/2}}\Biggl(\sum^{[t/(2v_n)]-1}_{j=0}
\bigl(\xi _{0,j}(u_{n})/u_{n}^2-E
\bigl(\xi_{0,j}(u_{n})/u^2_{n}|
\mathcal{F}_{2jv_n}\bigr)\bigr),
\nonumber
\\[-8pt]
\\[-8pt]
\nonumber
&&    \hspace*{10pt}\qquad \sum^{[t/(2v_n)]-1}_{j=0}\bigl(
\xi_{1,j}(u_{n})/u_{n}^2-E\bigl(\xi
_{1,j}(u_{n})/u_{n}^2|
\mathcal{F}_{2jv_n}\bigr)\bigr)\Biggr).
\end{eqnarray}
Let
\begin{eqnarray*}
\chi^{n, 0}_j&=&\frac{2v_n}{\Delta_n^{1/2}}\bigl(\xi
_{0,j}(u_{n})/u_{n}^2-E\bigl(
\xi_{0,j}(u_{n})/u^2_{n}|\mathcal
{F}_{2jv_n}\bigr)\bigr),
\\
\chi^{n,
1}_j&=&\frac{2v_n}{\Delta_n^{1/2}}\bigl(
\xi_{1,j}(u_{n})/u_{n}^2-E\bigl(\xi
_{1,j}(u_{n})/u^2_{n}|
\mathcal{F}_{2jv_n}\bigr)\bigr).
\end{eqnarray*}
By Theorem~7.28 in Chapter IX of Jacod and Shiyayev (\citeyear{JaSh03}), we only
need to prove that
%
\begin{equation}\quad
\label{stcon} \cases{ %
\displaystyle \sup_t\Biggl|\sum
^{[{t}/{(2v_n)}]-1}_{j=0}E\bigl(\chi^{n, k}_j|
\mathcal {F}_{2jv_n}\bigr)\biggr|\rightarrow^P 0; & \quad$ k=0, 1,$
\vspace*{2pt}\cr
\displaystyle \sum^{[{t}/{(2v_n)}]-1}_{j=0}E_{\mathcal{F}_{2jv_n}} \bigl(
\chi^{n,
k}_j \bigr)^2\rightarrow^P
4\int^t_0c_s^2\,ds; &\quad
$k=0, 1,$
\vspace*{2pt}\cr
\displaystyle \sum^{[{t}/{(2v_n)}]-1}_{j=0}E_{\mathcal{F}_{2jv_n}} \bigl(
\chi^{n,
0}_j\chi^{n, 1}_j \bigr)
\rightarrow^P 2\int^t_0c_s^2\,ds;
&
\vspace*{2pt}\cr
\displaystyle \sum^{[{t}/{(2v_n)}]-1}_{j=0}E_{\mathcal{F}_{2jv_n}} \bigl(
\chi^{n,
k}_j \bigr)^2I_{\{|\chi^{n, k}_j|>\varepsilon\}}
\rightarrow^P 0; &\quad  $k=0, 1,$
\vspace*{2pt}\cr
\displaystyle \sum^{[{t}/{(2v_n)}]-1}_{j=0}E_{\mathcal{F}_{2jv_n}} \bigl(
\chi^{n,
k}_j(M_{2(j+1)v_n}-M_{2jv_n}) \bigr)
\rightarrow^P 0; &\quad $k=0, 1,$}
\end{equation}
for any square-integrable martingale $M$. The first equation holds
automatically since $(\chi^{n, k}_j)_{j=0}^{[{t}/{(2v_n)}]-1}$ form a
sequence of $\mathcal{F}_{2(j+1)v_n}$-martingale differences.

Now we calculate the conditional variances of the marginal sequences.
By (\ref{cond}), Lemma~\ref{lem2} and the fact that
$|U_j(u_{n})-e^{-u_{n}^2\sigma^2_{(i-1)\Delta_n}}|\leq K\Delta
_n^{1-\beta/2}u_{n}^{\beta}$, we have
%
\begin{eqnarray}
\label{cv1} &&\frac{4v_n^2}{\Delta_nu_{n}^4}\sum^{[t/(2v_n)]-1}_{j=0}E
\bigl(\bigl(\xi _{k,j}(u_{n})-E\bigl(\xi_{k,j}(u_{n})|
\mathcal{F}_{2jv_n}\bigr)\bigr)^2|\mathcal {F}_{2jv_n}
\bigr)
\nonumber
\\
&&\qquad= \frac{4v_n^2}{\Delta_nu_{n}^4}\sum^{[t/(2v_n)]-1}_{j=0}
\bigl(E\bigl[\xi _{k,j}^2(u_{n})|
\mathcal{F}_{2jv_n}\bigr]-\bigl(E\bigl[\xi_{k,j}(u_{n})|
\mathcal {F}_{2jv_n}\bigr]\bigr)^2 \bigr)
\nonumber
\\
&&\qquad= \frac{4v_n^2}{\Delta_nu_{n}^4}\sum^{[t/(2v_n)]-1}_{j=0}E
\bigl(\xi _{k,j}^2(u_n)|\mathcal{F}_{2jv_n}
\bigr)+o_p(1)
\nonumber
\\
&&\qquad= \frac{4v_n^2}{2(k_n-1)\Delta_nu_{n}^4}\sum^{[t/(2v_n)]-1}_{j=0}
\frac
{U_j(2u_{n})+1-2U_j^2(u_{n})}{U_j^2(u_n)}+o_p(1)
\\
&&\qquad= \frac{(1+o_p(1))}{u^4_{n}}\sum^{[{t}/{(2v_n)}]-1}_{j=0}
\bigl(U_j(2u_{n})+1-2U_j^2(u_{n})
\bigr)2v_n+o_p(1)
\nonumber
\\
&&\qquad= \frac{(1+o_p(1))}{u^4_{n}}\sum^{[{t}/{(2v_n)}]-1}_{j=0}
\bigl(e^{-4u_{n}^2c_{2jv_n}}+1-2e^{-2u^2_{n}c_{2jv_n}} \bigr)2v_n+o_p(1)
\nonumber
\\
&&\qquad= \frac{\int^t_0(e^{-4u_{n}^2c_s}+1-2e^{-2u_{n}^2c_s})\,ds}{u_{n}^4}+o_p(1)\rightarrow ^P 4\int
^t_0c_s^2\,ds,\nonumber
\end{eqnarray}
where in obtaining the convergence in probability, we used the Taylor
expansion of $e^x$ when $x$ is near 0. This proves the second equation
in (\ref{stcon}).

Next, we are going to check the third equation in (\ref{stcon}). By
Lemma~\ref{lem2}, we have
%
\begin{eqnarray}
\label{covar} &&\frac{4v_n^2}{\Delta_nu_{n}^4}\sum^{[t/(2v_n)]-1}_{j=0}E_{\mathcal
{F}_{2jv_n}}
\bigl(\xi_{0,j}(u_{n})-E_{\mathcal{F}_{2jv_n}}\xi
_{0,j}(u_{n}) \bigr)
\nonumber
\\
&&\hspace*{59pt}\quad\times{} \bigl(\xi_{1,j}(u_{n})-E_{\mathcal{F}_{2jv_n}}\xi
_{1,j}(u_{n}) \bigr)
\nonumber
\\
&&\qquad= \frac{4v_n^2}{\Delta_nu_{n}^4}\sum^{[t/(2v_n)]-1}_{j=0}
\bigl(E_{\mathcal
{F}_{2jv_n}}\xi_{0,j}(u_{n})\xi_{1,j}(u_{n})\\
&&\hspace*{104pt}{}-E_{\mathcal{F}_{2jv_n}}
\xi _{0,j}(u_{n})E_{\mathcal{F}_{2jv_n}}\xi_{1,j}(u_{n})
\bigr)\nonumber
\\
&&\qquad=\frac{4v_n^2}{\Delta_nu_{n}^4}\sum^{[t/(2v_n)]-1}_{j=0}E_{\mathcal
{F}_{2jv_n}}
\xi_{0,j}(u_{n})\xi_{1,j}(u_{n})+o_p(1).\nonumber
\end{eqnarray}
Now we investigate the summand in (\ref{covar}). Let
\begin{eqnarray*}
\zeta_k(j, l)&=&\cos\biggl(u_{n}\frac{\Delta^n_{2jk_n+2l-k+1}X-\Delta
^n_{2jk_n+2l-k}X}{\Delta_n^{1/2}}
\biggr)
\\
&&{}-E_{\mathcal{F}_{(2jk_n+2l-k-1)\Delta_n}}\cos\biggl(u_{n}\frac{\Delta
^n_{2jk_n+2l-k+1}X-\Delta^n_{2jk_n+2l-k}X}{\Delta_n^{1/2}}\biggr)
\end{eqnarray*}
and
\[
\zeta_k^{\prime}(j, l)=\cos\biggl(u_{n}
\frac{\Delta^n_{2jk_n+2l-k+1}X-\Delta
^n_{2jk_n+2l-k}X}{\Delta_n^{1/2}}\biggr)-U_j(u_{n}),
\]
$k=0, 1$. By (6.22) and (6.29) in \citet{26}, we have
%
\begin{eqnarray}
\label{545}&& \bigl|\zeta^{\prime}_k(j, l)-\zeta_k(j,
l)\bigr|
\nonumber
\\[-8pt]
\\[-8pt]
\nonumber
&&\qquad\leq Ku^4_{n}\Delta_n^{1/2}
\phi_n+\bigl|U_{2jv_n+(2l-k-1)\Delta_n}(u_n)-U_j(u_{n})\bigr|,
\end{eqnarray}
which, together with Lemma~\ref{lem2} and the property of $U_t(u_n)$,
shows that\vspace*{-1pt}
%
\begin{eqnarray}
\label{ud1} & &\biggl |E_{\mathcal{F}_{2jv_n}}\xi_{1, j}(u_n) \biggl(
\frac{
({1}/{(k_n-1)})\sum^{k_n-1}_{l=1}(\zeta^{\prime}_0(j, l)-\zeta_0(j,
l))}{U_j(u_{n})} \biggr)\biggr|
\nonumber
\\[-1pt]
&&\qquad\leq \sqrt{E_{\mathcal{F}_{2jv_n}}\xi^2_{1, j}(u_n)}
\sqrt{E_{\mathcal
{F}_{2jv_n}} \biggl(\frac{({1}/{(k_n-1)})\sum^{k_n-1}_{l=1}\llvert \zeta
^{\prime}_0(j, l)-\zeta_0(j, l)\rrvert }{U_j(u_{n})} \biggr)^2}
\nonumber
\\[-8pt]
\\[-8pt]
\nonumber
&&\qquad \leq K \frac{u_{n}^2}{\sqrt{k_n}} \Bigl(u^4_{n}
\Delta_n^{1/2}\phi _n+\sqrt{\max
_{l}E_{\mathcal{F}_{2jv_n}}\bigl(U_{2jv_n+(2l-1)\Delta
_n}-U_j(u_{n})
\bigr)^2} \Bigr)
\\[-1pt]
&&\qquad \leq K\frac{u_n^4\sqrt{v_n}}{\sqrt{k_n}}.\nonumber\vspace*{-1pt}
\end{eqnarray}
Similarly, by the property of $U_t(u_n)$, (6.22) and (6.29)
in \citet{26}, and H\"{o}lder's inequality, we have\vspace*{-1pt}
%
\begin{eqnarray}
\label{ud2} &&\biggl|E_{\mathcal{F}_{2jv_n}}\frac{\sum^{k_n-1}_{l=1}\zeta_0(j,
l)}{(k_n-1)U_j(u_n)} \biggl(\frac{({1}/{(k_n-1)})\sum^{k_n-1}_{l=1}(\zeta
^{\prime}_1(j, l)-\zeta_1(j, l))}{U_j(u_{n})}
\biggr)\biggr|
\nonumber
\\[-1pt]
&&\qquad\leq \sqrt{\frac{\sum^{k_n-1}_{l=1}E_{\mathcal{F}_{2jv_n}}\zeta
_0^2(j, l)}{(k_n-1)^2U_j^2(u_n)}}
\nonumber
\\[-8pt]
\\[-8pt]
\nonumber
&&\qquad\quad{}\times \sqrt{E_{\mathcal{F}_{2jv_n}}
\biggl(\frac{({1}/{(k_n-1)})\sum^{k_n-1}_{l=1}\llvert \zeta^{\prime}_0(j,
l)-\zeta_0(j, l)\rrvert }{U_j(u_{n})} \biggr)^2}
\\[-1pt]
&&\qquad\leq K\frac{u_n^4\sqrt{v_n}}{\sqrt{k_n}}.\nonumber
\end{eqnarray}
Equations (\ref{ud1}) and (\ref{ud2}) yield\vspace*{-1pt}
%
\begin{eqnarray}\qquad
\label{covar1} &&E_{\mathcal{F}_{2jv_n}}\xi_{0,j}(u_{n})
\xi_{1,j}(u_{n})
\nonumber
\\[-8pt]
\\[-8pt]
\nonumber
&&\qquad=E_{\mathcal
{F}_{2jv_n}}\frac{({1}/{(k_n-1)})\sum^{k_n-1}_{l=1}\zeta_0(j, l)
({1}/{(k_n-1)})\sum^{k_n-1}_{l=1}\zeta_1(j, l)}{U_j(u_{n})U_j(u_{n})}+r_j,
\end{eqnarray}
where $r_j$ satisfies $|r_j|\leq K\sqrt{v_n}u^4_{n}/\sqrt{k_n}$. By the
definition of $\zeta_k(j, l)$, we have\vspace*{-1pt}
\begin{eqnarray}
\label{covar12} && E_{\mathcal{F}_{2jv_n}}\frac{({1}/{(k_n-1)})\sum^{k_n-1}_{l=1}\zeta
_0(j, l)({1}/{(k_n-1)})\sum^{k_n-1}_{l=1}\zeta_1(j,
l)}{U_j(u_{n})U_j(u_{n})}
\nonumber
\\[-1pt]
&&\qquad=\frac{1}{(k_n-1)^2 U_j(u_{n})U_j(u_{n})}\sum^{k_n-1}_{l=1}E_{\mathcal
{F}_{2jv_n}}
\zeta_0(j, l)\zeta_1(j, l)\\[-1pt]
&&\qquad\quad{}+\sum
^{k_n-2}_{l=1}E_{\mathcal
{F}_{2jv_n}}\zeta_0(j,l)
\zeta_1(j, l+1).\nonumber
\end{eqnarray}
By Lemmas 11--12 in \citet{26}, we have
%
\begin{eqnarray}
\label{covar2} && E_{\mathcal{F}_{2jv_n}}\zeta_0(j, l)\zeta_1(j,
l)
\nonumber
\\
&&\qquad= E_{\mathcal{F}_{2jv_n}}\zeta_0(j, l)\cos\biggl(u_{n}
\frac{\Delta
_{2jk_n+2l}^nX-\Delta^n_{2jk_n+2l-1}X}{\Delta_n^{1/2}}\biggr)
\nonumber
\\
&&\qquad=E_{\mathcal{F}_{2jv_n}}\cos\biggl(u_{n}\frac{\Delta_{2jk_n+2l+1}^nX-\Delta
^n_{2jk_n+2l}X}{\Delta_n^{1/2}}\biggr)\\
&&\qquad\quad{}\times \cos
\biggl(u_{n}\frac{\Delta
_{2jk_n+2l}^nX-\Delta^n_{2jk_n+2l-1}X}{\Delta_n^{1/2}}\biggr)
\nonumber
\\
&&\qquad\quad{}-U_{j}(u_{n})U_{j}(u_{n})+r_{2j},\nonumber
\end{eqnarray}
where $r_{2j}$ satisfies $|r_{2j}|\leq Ku_n^2\sqrt{v_n}$. Since
$\cos(x)\cos(y)=\frac{1}{2}(\cos(x+y)+\cos(x-y))$, we have by Lemma~\ref{lem3},
%
\begin{eqnarray}
\label{covar3} &&\hspace*{-4pt}E_{\mathcal{F}_{2jv_n}}\cos\biggl(u_{n}\frac{\Delta_{2jk_n+2l+1}^nX-\Delta
^n_{2jk_n+2l}X}{\Delta_n^{1/2}}
\biggr)\cos\biggl(u_{n}\frac{\Delta
_{2jk_n+2l}^nX-\Delta^n_{2jk_n+2l-1}X}{\Delta_n^{1/2}}\biggr)\hspace*{-60pt}
\nonumber
\\
&&\hspace*{-6pt}\qquad= \frac{1}{2} E_{\mathcal{F}_{2jv_n}}\biggl(\cos\biggl(u_{n}
\frac{\Delta
^n_{2jk_n+2l+1}X}{\Delta_n^{1/2}}-u_{n}\frac{\Delta
^n_{2jk_n+2l-1}X}{\Delta_n^{1/2}}\biggr)
\nonumber
\\
&&\hspace*{-6pt}\hspace*{38pt}\qquad\quad{} + \cos\biggl(u_{n}\frac{\Delta^n_{2jk_n+2l+1}X}{\Delta
_n^{1/2}}-2u_{n}
\frac{\Delta^n_{2jk_n+2l}X}{\Delta_n^{1/2}}+u_{n}\frac
{\Delta^n_{2jk_n+2l-1}X}{\Delta_n^{1/2}}\biggr) \biggr)
\\
&&\hspace*{-6pt}\qquad= \frac{1}{2}E_{\mathcal{F}_{2jv_n}} \bigl(\exp {\bigl(-u_n^2c_{2jv_n+(2l-2)\Delta_n}
\bigr)}+\exp{\bigl(-3u_n^2c_{2jv_n+(2l-2)\Delta
_n}\bigr)}
\bigr)+r_{3, j}
\nonumber
\\
&&\hspace*{-6pt}\qquad= \frac{1}{2} \bigl(\exp{\bigl(-u_n^2c_{2jv_n}
\bigr)}+\exp {\bigl(-3u_n^2c_{2jv_n}\bigr)}
\bigr)+r_{3, j}+r_{4, j},\nonumber
\end{eqnarray}
where $|r_{3,j}|\leq K \Delta_n^{1-\beta/2}$ and $|r_{4, j}|\leq
Ku_n^2v_n$ by second-order Taylor expansion on $e^x$ for $x$ around the
origin and (S.1.3) with $V=c$. Now substituting (\ref{covar3}) back
into (\ref{covar2}), we have
%
\begin{eqnarray}
\label{covar4} E_{\mathcal{F}_{2jv_n}}\zeta_0(j, l)\zeta_1(j,
l)&=&\tfrac{1}{2} \bigl(\exp{\bigl(-u_n^2c_{2jv_n}
\bigr)}+\exp{\bigl(-3u_n^2c_{2jv_n}\bigr)} \bigr)
\nonumber
\\[-8pt]
\\[-8pt]
\nonumber
&&{} -\exp{\bigl(-2u_n^2c_{2jv_n}
\bigr)}+r_{5, j},
\end{eqnarray}
where $|r_{5j}|\leq K (\sqrt{v_n}+\Delta_n^{1-\beta/2})$.
Similarly, we have
%
\begin{eqnarray}
\label{covar5} && E_{\mathcal{F}_{2jv_n}}\zeta_0(j, l)\zeta_1(j,
l+1)\nonumber\\
&&\qquad=\tfrac{1}{2} \bigl(\exp{\bigl(-u_n^2c_{2jv_n}
\bigr)}+\exp{\bigl(-3u_n^2c_{2jv_n}\bigr)} \bigr)
\\
&&\qquad\quad{} -\exp{\bigl(-2u_n^2c_{2jv_n}
\bigr)}+r_{6, j},\nonumber
\end{eqnarray}
where $|r_{6,j}|\leq K (\sqrt{v_n}+\Delta_n^{1-\beta/2}) $. Substitute
(\ref{covar4}) and (\ref{covar5}) into (\ref{covar12}), and then
substitute the latter into (\ref{covar1}), and we have
%
\begin{eqnarray}
\label{covar6} & & E_{\mathcal{F}_{2jv_n}}\xi_{0,j}(u_{n})
\xi_{1,j}(u_{n})
\nonumber
\\
&&\qquad= \frac{\exp{(-u_n^2c_{2jv_n})} +\exp{(-3u_n^2c_{2jv_n})}-2\exp
{(-2u_n^2c_{2jv_n})}}{(k_n-1)U_j^2(u_{n})}+r^*_{j}
\\
&&\qquad= \frac{u_n^4c_{2jv_n}^2+r_{7j}}{(k_n-1)U_j^2(u_{n})}+r^*_{j}=\frac
{u_n^4c_{2jv_n}^2+r_{8j}}{k_n-1}+r^*_{j},\nonumber
\end{eqnarray}
where $|r_{7j}|\vee|r_{8j}|\leq Ku_n^6$, $|r^*_j|\leq|r_j|+\frac
{|r_{2j}|+|r_{3j}|+|r_{4j}|+|r_{5j}|+|r_{6j}|}{k_n-1}+\frac
{1}{(k_n-1)^2}$. Now a combination of (\ref{covar6}) and (\ref{covar}) yields
%
\begin{eqnarray}
\label{covar7} &&\frac{4v_n^2}{\Delta_nu_n^4}\sum^{[t/(2v_n)]-1}_{j=0}E_{\mathcal
{F}_{2jv_n}}
\bigl(\xi_{0,j}(u_{n})-E_{\mathcal{F}_{2jv_n}}\xi
_{0,j}(u_{n}) \bigr)
\nonumber
\\[-8pt]
\\[-8pt]
\nonumber
&&\hspace*{48pt}\qquad{}\times \bigl(\xi_{1,j}(u_{n})-E_{\mathcal{F}_{2jv_n}}\xi
_{1,j}(u_{n}) \bigr)\rightarrow^p 2\int
^t_0c_s^2\,ds.
\end{eqnarray}
This proves the third equation in (\ref{stcon}).

By Lemma~\ref{lem2}, we also have
\begin{eqnarray*}
&&\sum^{[{t}/{(2v_n)}]-1}_{j=0}E_{\mathcal{F}_{2jv_n}}\bigl(
\chi^{n,
k}_j\bigr)^2I\bigl(\bigl|
\chi^{n, k}_j\bigr|>\varepsilon\bigr)\\
&&\qquad\leq\frac{1}{\varepsilon}\sum
^{[
{t}/{(2v_n)}]-1}_{j=0}E_{\mathcal{F}_{2jv_n}}\bigl|
\chi^{n, k}_j\bigr|^3
\\
&&\qquad\leq \frac{K}{\varepsilon}\sum_{j=0}^{[{t}/{(2v_n)}]-1}
\biggl(\frac{2v_n}{\Delta
_n^{1/2}}\biggr)^3\frac{1}{u_{n}^6}E_{\mathcal{F}_{2jv_n}}\bigl|
\xi_{k,
j}(u_{n})\bigr|^3 \rightarrow0.
\end{eqnarray*}
This proves the Linderberg condition [equation four in (\ref{stcon})].

Taking $\kappa=2$ and $\zeta^n_j=1$ in Lemma~15 of \citet{26}, we have
\begin{eqnarray*}
&&\sum^{[{t}/{(2v_n)}]-1}_{j=0}E_{\mathcal{F}_{2jv_n}}
\chi^{n,
k}_j(M_{2(j+1)v_n}-M_{2jv_n})
\\
&&\qquad=\sum^{[{t}/{(2v_n)}]-1}_{j=0}\frac{2v_n}{u_{n}^2\Delta
_n^{1/2}}E_{\mathcal{F}_{2jv_n}}
\xi_{k,
j}(u_{j}) (M_{2(j+1)v_n}-M_{2jv_n})
\rightarrow^P 0.
\end{eqnarray*}
This proves the final equation in (\ref{stcon}) and completes the proof
of the bivariate central limit theorem with stable convergence.
\end{pf*}

\begin{pf*}{Proof of Theorem~\ref{th2}} Let $T_{n1}=\hat
{C}_0(u_n)-A_0( u_n)^n_t-C_t-(\hat{C}_1(u_n)-A_0(u_n)^n_t-C_t)$. By
(\ref{decom}),
%
\begin{equation}
\label{statist} T_n\equiv\frac{T_{n1}-\gamma_n\Delta
_n^{1/2}}{C_t+A_0(u_n)^n_t+O_p(\Delta_n^{1/2})}=\frac
{T_{n1}}{C_t+o_p(1)}+o_p
\bigl(\Delta_n^{1/2}\bigr).
\end{equation}
Then Theorem~\ref{th2} is a straightforward consequence of Theorem~\ref
{th1}, (\ref{statist}), the stable convergence mode and the continuous
mapping theorem.
\end{pf*}

\begin{pf*}{Proof of Theorem~\ref{th3}} By Theorem~\ref{th1},
$\hat{C}_1(u_n)=C_t+A_0(u_n)^n_t+O_p(\Delta_n^{1/2})=C_t+o_p(1)$. This shows
that the denominator of $\hat{\kappa}_T$ converges to $C_T^2$ in
probability. By (\ref{d1}), we have
%
\begin{eqnarray}
&&\biggl(c^k_j(u_n)-\frac{(\sinh
(u_n^2c^k_j(u_n)))^2}{(k_n-1)u_n^2}
\biggr)^2
\nonumber
\\[-8pt]
\\[-8pt]
\nonumber
&&\qquad=c^2_{2jv_n}+\tilde{c}^k_{j,
1}(u_n)+
\tilde{c}^k_{j, 2}(u_n)+\tilde{c}^k_{j, 3}(u_n),
\end{eqnarray}
where
\begin{eqnarray*}
\tilde{c}^k_{j, 1}&=&\biggl(\frac{\xi_{k,
j}(u_n)}{u_n^2}
\biggr)^2+\biggl(\frac{\xi^2_{k,
j}(u_n)}{2u_n^2}-\frac{(\sinh(u_n^2c^k_j(u_n)))^2}{(k_n-1)u_n^2}
\biggr)^2+\bigl(r_{k,
j}(u_n)
\bigr)^2,
\\
\tilde{c}^k_{j, 2}&=&2c_{2jv_n}\biggl(-
\frac{\xi_{k,
j}(u_n)}{u_n^2}+\frac{\xi^2_{k,
j}(u_n)}{2u_n^2}-\frac{(\sinh(u_n^2c^k_j(u_n)))^2}{(k_n-1)u_n^2}+r_{k,
j}(u_n)
\biggr),
\\
\tilde{c}^k_{j, 3}(u_n)&=& 4c_{2jv_n}u_n^{\beta-2}
\Delta_n^{1-\beta
/2}a_{2jv_n}+\bigl(2u_n^{\beta-2}
\Delta_n^{1-\beta/2}a_{2jv_n}\bigr)^2.
\end{eqnarray*}
By (\ref{omega}),
\[
\sum^{[{t}/{(2v_n)}]-1}_{j=0}2v_n
\tilde{c}^k_{j, 1}I_{\Omega^c(k, n,
t)}=o_p(1),\qquad
\sum^{[{t}/{(2v_n)}]-1}_{j=0}2v_n
\tilde{c}^k_{j,
2}I_{\Omega^c(k, n, t)}=o_p(1).
\]
By Lemma~\ref{lem2},
$\sum^{[{t}/{(2v_n)}]-1}_{j=0}2v_n(\frac{\xi_{k,
j}(u_n)}{u_n^2})^2=o_p(1)$. On $\Omega(k, n, t)$, $|\frac{\xi^2_{k,
j}(u_n)}{2u_n^2}-\frac{(\sinh(u_n^2c^k_j(u_n)))^2}{(k_n-1)u_n^2}|$
is bounded by $K/u_n^2$, hence
\begin{eqnarray}
&&\biggl(\frac{\xi^2_{k, j}(u_n)}{2u_n^2}-\frac{(\sinh
(u_n^2c^k_j(u_n)))^2}{(k_n-1)u_n^2}\biggr)^2I_{\Omega(k, n, t)}\nonumber
\\
&&\qquad\leq \frac
{K}{u_n^2}\biggl\llvert \frac{\xi^2_{k, j}(u_n)}{2u_n^2}-\frac{(\sinh
(u_n^2c^k_j(u_n)))^2}{(k_n-1)u_n^2}\biggr
\rrvert I_{\Omega(k, n, t)}
\nonumber
\\
&& \qquad\leq\frac{K}{u_n^2}\biggl(\biggl|\frac{\xi^2_{k,
j}(u_n)}{2u_n^2}-E_{\mathcal{F}_{2jv_n}}
\frac{\xi^2_{k,
j}(u_n)}{2u_n^2}\biggr|\\
&&\hspace*{20pt}\qquad\quad{}+\biggl|E_{\mathcal{F}_{2jv_n}}\frac{\xi^2_{k,
j}(u_n)}{2u_n^2}-\frac
{U_j(2u_n)+1-2U_j^2(u_n^2)}{4u_n^2(k_n-1)U_j^2(u_n)}\biggr|
\nonumber
\\
& &\hspace*{20pt}\qquad\quad{} +\biggl|\frac{U_j(2u_n)+1-2U_j^2(u_n^2)}{4u_n^2(k_n-1)U_j^2(u_n)}-\frac{(\sinh
(u_n^2c^k_j(u_n)))^2}{(k_n-1)u_n^2}\biggr|\biggr)I_{\Omega(k,
n, t)}.
\nonumber
\end{eqnarray}
By the property of $U_j(u_n)$ and the definition of $c^k_j(u_n)$, the
expectation of the third absolute value conditional
on $\mathcal{F}_{2jv_n}$ is smaller than
$K(u_n^{\beta-2}\Delta_n^{1-\beta/2}/k_n+u_n^{-2}/k_n^{3/2}+\Delta
_n^{1/2}\phi_n/k_n)$.
By Lemma~\ref{lem2}, the second absolute value is smaller than
$Ku_n^4\Delta_n^{1/2}\phi_n$. By H\"{o}lder's inequality and Lemma~\ref{lem2}
with $q=4$, the expectation of the first absolute value conditional
on $\mathcal{F}_{2jv_n}$ is smaller than $K(u_n^2/k_n+\sqrt{v_n})$.
In summary, we conclude that
%
\begin{equation}
\sum^{[{t}/{(2v_n)}]-1}_{j=0}2v_n\biggl(
\frac{\xi^2_{k,
j}(u_n)}{2u_n^2}-\frac{(\sinh(u_n^2c^k_j(u_n)))^2}{k_nu_n^2}\biggr)^2I_{\Omega(k,
n, t)}=o_p(1).
\end{equation}
By (\ref{omega}), $\sum_{j=0}^{[{t}/{(2v_n)}]-1}2v_n(r_{k,
j}(u_n))^2I_{\Omega^c(k, n, t)}=o_p(1)$. On $\Omega(k, n, t)$, $|r_{k,
j}|\leq K\frac{|\xi_{k, j}(u_n)|^3}{u_n^2}$. By Lemma~\ref{lem2} with
$q=6$, we have $\sum_{j=0}^{[{t}/{(2v_n)}]-1}2v_n(r_{k,
j}(u_n))^2\times\break  I_{\Omega(k, n, t)}=o_p(1)$. Combining all the results of the
terms on the right-hand side of the decomposition of $\tilde{c}^k_{j,
1}$, we have $\sum_{j=0}^{[{t}/{(2v_n)}]-1}2v_n\tilde{c}^k_{j,
1}=o_p(1)$. Similarly, one easily proves that $\sum_{j=0}^{[
{t}/{(2v_n)}]-1}2v_n\tilde{c}^k_{j, 2}=o_p(1)$. By boundedness of $c$ and
$a$, $\sum^{[{t}/{(2v_n)}]-1}_{j=0}2v_n\tilde{c}^k_{j, 3}=o_p(1)$.
This shows that
%
\begin{equation}\qquad
\label{asyvar} \hat{I}_{nk}=\sum^{[{T}/{(2v_n)}]-1}_{j=0}2v_nc_{2jv_n}^2+o_p(1)=
\int^T_0c_s^2\,ds+o_p(1),\qquad
k=0, 1.
\end{equation}
This shows that the numerator of $\hat{\kappa}_T$ converges to
$4\int^T_0c_s^2\,ds$ in probability, and hence $\hat{\kappa}_T$ itself
converges to $\kappa_T$ in probability. On the other hand, by
Theorem~\ref{th2}, $T_n/\Delta_n^{1/2}$ converges to $G_T$ stably.
By the stable convergence mode, $\mathcal{T}_n$ converges to
standard normal distribution stably.
\end{pf*}

\subsection{Proof of results under $H_1$}
In the sequel we assume that $X$ is a pure-jump process. We rewrite
%
\begin{eqnarray}
\label{decom1} c^k_j(u_n)&=&-
\frac{\log{U_j(u_n)}}{u_n^2}-\frac{\log{(1+\xi
_{k,j}(u_n))}}{u_n^2}
\nonumber
\\[-8pt]
\\[-8pt]
\nonumber
&=&2u_n^{\beta-2}\Delta_n^{1-\beta/2}a_{2jv_n}-
\frac{\xi
_{k,j}(u_n)}{u_n^2}+\tilde{r}_{k, j}, \qquad k=0, 1,
\end{eqnarray}
where $|\tilde{r}_{k,j}|\leq K \xi^2_{k,j}(u_n)/u_n^2$ on $\Omega(k, n,
t)$. Recall the definition of $T_{n1}$ in (\ref{statist}), and we have
%
\begin{equation}
\label{diff} T_{n1}=2v_n\sum
^{[t/(2v_n)]-1}_{j=0}-\frac{\xi_{0,j}(u_n)-\xi
_{1,j}(u_n)}{u_n^2}+\tilde{R}_{n,t},
\end{equation}
where
%
\begin{eqnarray}
\label{remain} &&\tilde{R}_{n,t}=2v_n\sum
^{[t/(2v_n)]-1}_{j=0}\biggl[(\tilde{r}_{0,
j}-
\tilde{r}_{1,
j})
\nonumber
\\[-8pt]
\\[-8pt]
\nonumber
&&\hspace*{97pt}{}+\biggl(\frac{(\sinh(u_n^2c^1_j(u_n)))^2-(\sinh
(u_n^2c^0_j(u_n)))^2}{u^2_n(k_n-1)}\biggr)\biggr].
\end{eqnarray}

Similar to Lemma~\ref{lem2}, we have the following. The proof is
provided in the supplementary material [\citet{KonLiuJin14}].

\begin{lemma}\label{lem40}
Assume Assumptions \ref{ass1}, \ref{ass3} and
\ref{ass4}, and suppose $u_n$ is bounded, so
we have on the set $\{C_t=0\}$,
%
\begin{eqnarray}
\label{h1} &&\bigl|E_{\mathcal{F}_{2jv_n}}\xi_{k, j}(u_n)\bigr|
\nonumber
\\
&&\qquad\leq K \bigl(\Delta_n^{(1-r/2)\wedge({(3-\beta)}/{2}-\varepsilon
^{\prime}) \wedge(1-{\beta}/{(2(\beta+1-r)}))}\\
&&\hspace*{72pt}{}+u_n^{\beta}
\Delta _n^{1-\beta/2}v_n^{\beta/2}+u_n^{2\beta}
\Delta_n^{2-\beta}v_n \bigr),\nonumber
\end{eqnarray}
and if further $k_n\Delta_n^{1/2-\varepsilon}\rightarrow\infty$ for any
$\varepsilon>0$, and $\sup_n\frac{k_n\Delta_n^{1/2}}{u_n^{4}}<\infty$ is
satisfied,
%
\begin{equation}
\label{h1square} E_{\mathcal{F}_{2jv_n}}\xi_{k, j}^2(u_n)
\leq K\frac{u_n^{\beta}\Delta
_n^{1-\beta/2}}{k_n}.
\end{equation}
\end{lemma}

The following lemma gives the convergence rate of the terms on the
right-hand side of (\ref{diff}). The proof can be found in the
supplementary material [\citet{KonLiuJin14}] to this paper.

\begin{lemma}\label{lem4} Assume Assumption~\ref{ass1}, \ref{ass3} and
\ref{ass4}, and suppose $u_n$ is bounded and $k_n\Delta
_n^{1/2}\rightarrow0$, so we have on the set $\{C_t=0\}$,
\begin{longlist}[(3)]
\item[(1)]
%
\begin{eqnarray}\qquad
&& \Biggl|2v_n\sum^{[t/(2v_n)]-1}_{j=0}E_{\mathcal{F}_{2jv_n}}
\biggl(\frac{\xi
_{0,j}(u_n)-\xi_{1,j}(u_n)}{u_n^2}\biggr)\Biggr|
\nonumber
\\[-8pt]
\\[-8pt]
\nonumber
&& \qquad\leq Ku_n^{-2}\bigl(\Delta_n^{(1-r/2)\wedge({(3-\beta)}/{2}-\varepsilon
^{\prime})\wedge(1-{\beta}/{(2(\beta+1-r))})}+u_n^{\beta}
\Delta _n^{3/2-\beta/2}\bigr);
\end{eqnarray}
\item[(2)]
%
\begin{eqnarray}\qquad
&&\Biggl|\sum^{[t/(2v_n)]-1}_{j=0}E_{\mathcal{F}_{2jv_n}}\biggl(
\frac{\xi
_{0,j}(u_n)-\xi_{1,j}(u_n)-E_{\mathcal{F}_{2jv_n}}(\xi_{0,j}(u_n)-\xi
_{1,j}(u_n))}{u_n^2/(2v_n)}\biggr)^2\Biggr|
\nonumber
\\[-8pt]
\\[-8pt]
\nonumber
&&\qquad\leq Ku_n^{-4}\bigl(u_n^{\beta
}
\Delta_n^{2-\beta/2}+\Delta_n^{(2-r/2)\wedge({(5-\beta)}/{2}-\varepsilon^{\prime})}\bigr);
\nonumber
\end{eqnarray}
for any $\varepsilon^{\prime}>0$;
\item[(3)]
%
\begin{equation}
\tilde{R}_{nt}=O_p\biggl(\frac{u_n^{\beta-2}\Delta_n^{1-\beta/2}}{k_n}\biggr).
\end{equation}
\end{longlist}
\end{lemma}

\begin{pf*}{Proof of Theorem~\ref{th4}} We first prove the first
equation. By (\ref{decom1}), we have
%
\begin{eqnarray}
\hat{C}_0(u_n)-\hat{C}_1(u_n)&=&T_{n, 1}
\nonumber
\\[-8pt]
\\[-8pt]
\nonumber
&=&2v_n
\sum^{[
{t}/{(2v_n)}]-1}_{j=0}-\frac{\xi_{0, j}(u_n)-\xi_{1, j}(u_n)}{u_n^2}+\tilde
{R}_{n, t}.
\end{eqnarray}
Now by Lemma~\ref{lem4}, we have
%
\begin{equation}
\label{fbias} \hat{C}_0(u_n)-\hat{C}_1(u_n)=O_p(
\delta_{n, 1}+\sqrt{\delta_{n,
2}}+\delta_{n, 3}),
\end{equation}
where
\begin{eqnarray*}
\delta_{n, 1}&=&u_n^{-2}\bigl(
\Delta_n^{(1-r/2)\wedge({(3-\beta)}/{2}-\varepsilon^{\prime})\wedge(1-{\beta}/{(2(\beta+1-r))})}+u_n^{\beta
}
\Delta_n^{3/2-\beta/2}\bigr),
\\
\sqrt{\delta_{n, 2}}&=&u_n^{-2}
\bigl(u_n^{\beta/2}\Delta_n^{1-\beta/4}+\Delta
_n^{(1-r/4)\wedge({(5-\beta)}/{4}-\varepsilon^{\prime}/2)}\bigr),\\
 \delta _{n, 3}&=&u_n^{\beta-2}
\Delta_n^{1-\beta/2}/k_n.
\end{eqnarray*}
Now, notice that: (1) $1-r/4>1-r/2>1-\beta/2(\beta+1-r)$; (2)
$\delta_{n, 3}>u_n^{\beta-2}\Delta_n^{3/2-\beta/2}$; (3)
$\frac{3-\beta}{2}<\frac{5-\beta}{4}$; (4) $\frac{\delta_{n,
3}}{u_n^{-2}\Delta_n^{3/2-\beta/2-\varepsilon^{\prime}}}\leq
u_n^{\beta}/(k_n\Delta_n^{1/2-\varepsilon^{\prime}})\leq K$;\break (5)~$\frac{u_n^{-2}\Delta_n^{3/2-\beta/2-\varepsilon^{\prime}}}{u_n^{\beta
/2-2}\Delta_n^{1-\beta/4}}=u_n^{4-\beta/2}\frac{k_n\Delta
_n^{1/2}}{u_n^4}\frac{1}{k_n\Delta_n^{\beta/4+\varepsilon^{\prime}}}\leq
K$. By choosing $\varepsilon^{\prime}>0$ small enough and the conditions
on $u_n$ and $k_n$, we have
\[
\hat{C}_0(u_n)-\hat{C}_1(u_n)=O_p
\bigl(u_n^{-2}\Delta_n^{1-{\beta
}/{(2(\beta+1-r))}}+
u_n^{\beta/2-2}\Delta_n^{1-\beta/4}\bigr).
\]

Next, we prove the second equation. By (\ref{decom1}), we have
%
\begin{eqnarray}
& & \bigl(c^k_j(u_n)\bigr)^2-
\bigl(2u_n^{\beta-2}\Delta_n^{1-\beta
/2}a_{2jv_n}
\bigr)^2-\biggl(\frac{\xi_{k, j}(u_n)}{u_n^2}\biggr)^2-(
\tilde{r}_{k,
j})^2
\nonumber
\\[-8pt]
\\[-8pt]
\nonumber
&&\qquad=  2\bigl(2u_n^{\beta-2}\Delta_n^{1-\beta/2}a_{2jv_n}
\bigr) \biggl(-\frac{\xi_{k,
j}(u_n)}{u_n^2}+\tilde{r}_{k, j}\biggr)-2\biggl(
\frac{\xi_{k, j}(u_n)}{u_n^2}\biggr)\tilde {r}_{k, j}.
\end{eqnarray}
Now we use several steps to show that under $H_1$ the principal term of
$(c^k_j(u_n))^2$ is $(2u_n^{\beta-2}\Delta_n^{1-\beta/2}a_{2jv_n})^2$
and $E_{\mathcal{F}_{2jv_n}}(c^k_j(u_n))^2\leq Ku_n^{2\beta-4}\Delta
_n^{2-\beta}$. By Lemma~\ref{lem40}, we have
%
\begin{equation}
\label{r1} E_{\mathcal{F}_{2jv_n}}\biggl(\frac{\xi_{k, j}(u_n)}{u_n^2}\biggr)^2
\leq K \bigl(u_n^{\beta-2}\Delta_n^{1-\beta/2}
\bigr)^2\frac{k_n\Delta
_n^{1/2}}{u_n^4}\frac{u_n^{4-\beta}}{k_n^2\Delta_n^{3/2-\beta/2}},
\end{equation}
which is $o_p((u_n^{\beta-2}\Delta_n^{1-\beta/2})^2)$ by the conditions
on $k_n$ and $u_n$ given in Theorem~\ref{th4}. By Lemma~\ref{lem40} and
(\ref{r1}), we have on $\Omega(k, n, t)$ (on which $|\tilde{r}_{k,
j}|\leq K\xi_{k, j}^2/u_n^2$ and $|\xi_{k, j}|$ is bounded),
%
\begin{eqnarray}
\label{r2} E_{\mathcal{F}_{2jv_n}}(\tilde{r}_{k, j})^2I_{\Omega(k, n, t)}
&\leq& KE_{\mathcal{F}_{2jv_n}}|\tilde{r}_{k, j}|I_{\Omega(k,n,t)}
\nonumber
\\
&\leq& KE_{\mathcal{F}_{2jv_n}}\biggl(\frac{\xi_{k,j}(u_n)}{u_n^2}\biggr)^2\\
&\leq& K
\bigl(u_n^{\beta-2}\Delta_n^{1-\beta/2}
\bigr)^2\frac{k_n\Delta
_n^{1/2}}{u_n^4}\frac{u_n^{4-\beta}}{k_n^2\Delta_n^{3/2-\beta/2}}.\nonumber
\end{eqnarray}
By (\ref{r1}) and (\ref{r2}), we have by H\"{o}lder's inequality,
%
\begin{eqnarray}
\label{r3} &&E_{\mathcal{F}_{2jv_n}}\biggl|\bigl(2u_n^{\beta-2}
\Delta_n^{1-\beta
/2}a_{2jv_n}\bigr) \biggl(-
\frac{\xi_{k, j}(u_n)}{u_n^2}+\tilde{r}_{k, j}\biggr)\biggr|I_{\Omega
(k, n, t)}
\nonumber
\\[-8pt]
\\[-8pt]
\nonumber
&&\qquad\leq K\bigl(u_n^{\beta-2}\Delta_n^{1-\beta/2}
\bigr)^2\biggl(\frac{u_n^{4-\beta
/2}}{k_n^{3/2}\Delta_n^{1-\beta/4}}+\frac{u_n^2}{k_n^2\Delta_n^{1/2}}\biggr)
\end{eqnarray}
and
%
\begin{eqnarray}
\label{r4} E_{\mathcal{F}_{2jv_n}}\biggl|\frac{\xi_{k, j}(u_n)}{u_n^2}\biggr|
 \tilde{r}_{k,
j}|I_{\Omega(k, n, t)}
&\leq& KE_{\mathcal{F}_{2jv_n}}\biggl(\frac{\xi
_{k,j}(u_n)}{u_n^2}\biggr)^2
\nonumber
\\[-8pt]
\\[-8pt]
\nonumber
&\leq& K \bigl(u_n^{\beta-2}\Delta_n^{1-\beta/2}
\bigr)^2\frac{k_n\Delta
_n^{1/2}}{u_n^4}\frac{u_n^{4-\beta}}{k_n^2\Delta_n^{3/2-\beta/2}}.
\end{eqnarray}
Combining (\ref{r1})--(\ref{r4}) yields that
%
\begin{equation}
\label{csquare} E_{\mathcal{F}_{2jv_n}}\bigl|\bigl(c^k_j(u_ n)
\bigr)^2-\bigl(2u_n^{\beta-2}\Delta_n^{1-\beta
/2}
\bigr)^2\bigr|I_{\Omega(k, n, t)}=o(1),
\end{equation}
where $o(1)$ holds uniformly in $j$.

By the form of $c^k_j(u_n)$, we have $u_n^2|c^k_j(u_n)|I_{\Omega(k, n,
t)}\leq K$, and hence by Taylor expansion on the exponential function,
we have
%
\begin{equation}
\label{sinhb} \bigl(\sinh\bigl(u_n^2c^k_j(u_n)
\bigr)\bigr)^2I_{\Omega(k, n, t)}\leq Ku_n^4
\bigl(c^k_j(u_n)\bigr)^2I_{\Omega(k, n, t)}
\leq K.
\end{equation}
By virtue of (\ref{sinhb}), we have
\begin{eqnarray}
\label{rsinh} &&E_{\mathcal{F}_{2jv_n}} \biggl(\frac{(\sinh(u_n^2c_j^k(u_n)))^2I_{\Omega
(k, n, t)}}{u_n^2(k_n-1)}
\biggr)^2\nonumber\\
&&\qquad\leq\frac{K}{u_n^4k_n^2}E_{\mathcal
{F}_{2jv_n}}\bigl(\sinh
\bigl(u_n^2c_j^k(u_n)
\bigr)\bigr)^2I_{\Omega(k, n, t)}
\\
&&\qquad\leq \frac{K}{k_n^2}E_{\mathcal{F}_{2jv_n}}\bigl(c_j^k(u_n)
\bigr)^2I_{\Omega
(k, n, t)}\leq\frac{K}{k_n^2}\bigl(u_n^{\beta-2}
\Delta_n^{1-\beta
/2}\bigr)^2,
\nonumber
\end{eqnarray}
and further by the Cauchy inequality,
%
\begin{equation}
\label{csinh} E_{\mathcal{F}_{2jv_n}}\biggl|c^k_j(u_n)
\frac{(\sinh
(u_n^2c_j^k(u_n)))^2}{u_n^2(k_n-1)}\biggr|I_{\Omega(k, n, t)}\leq \frac{K}{k_n}\bigl(u_n^{\beta-2}
\Delta_n^{1-\beta/2}\bigr)^2.
\end{equation}
Now combining (\ref{csquare}), (\ref{rsinh}), (\ref{csinh}) and (\ref
{omega}), we have
\begin{eqnarray*}
\hat{I}_{n, k}&=&\bigl(2u_n^{\beta-2}
\Delta_n^{1-\beta/2}\bigr)^2\Biggl(\sum
^{[
{t}/{(2v_n)}]-1}_{j=0}2v_na_{2jv_n}^2+o_p(1)
\Biggr)\\
&=&\bigl(2u_n^{\beta-2}\Delta _n^{1-\beta/2}
\bigr)^2\biggl(\int^t_0a_s^2\,ds+o_p(1)
\biggr),
\end{eqnarray*}
for $k=0, 1$. This proves the second equation of Theorem~\ref{th4}.
\end{pf*}

\begin{pf*}{Proof of Corollary~\ref{cor1}} Part 1 is a straight
consequence of Theorem~\ref{th3}. Now we prove part 2. By Theorem~\ref
{th4}, we have by the condition on $\gamma_n$,
%
\begin{eqnarray}
\mathcal{T}_{n}&=&\frac{-\gamma_n+O_p(u_n^{-2}\Delta_n^{1/2-{\beta
}/{(2(\beta+1-r))}}+ u_n^{\beta/2-2}\Delta_n^{1/2-\beta/4})}{4u_n^{\beta
-2}\Delta_n^{1-\beta/2}\sqrt{\int^t_0a_s^2\,ds+o_p(1)}}
\nonumber
\\[-8pt]
\\[-8pt]
\nonumber
&=&\frac{-\gamma_n(1+o_p(1))}{4u_n^{\beta-2}\Delta_n^{1-\beta/2}\sqrt
{\int^t_0a_s^2\,ds+o_p(1)}}.
\end{eqnarray}
Since $\gamma_nu_n^{2-\beta/2}\Delta_n^{{\beta}/{4}-
{1}/{2}}\rightarrow\infty$ and $\frac{u_n^{2-\beta/2}\Delta_n^{
{\beta}/{4}-{1}/{2}}}{u_n^{2-\beta}\Delta_n^{\beta/2-1}}\leq
u_n^{\beta/2}\Delta_n^{1/2-\beta/4}\rightarrow0$,
\[
\frac{-\gamma_n(1+o_p(1))}{4u_n^{\beta-2}\Delta_n^{1-\beta/2}\sqrt{\int^t_0a_s^2\,ds+o_p(1)}}\rightarrow^P -\infty.
\]
This proves part 2 on the performance of the power of the test.
\end{pf*}
\end{appendix}

\begin{supplement}[id=suppA]
\stitle{Supplement to ``Testing for pure-jump processes for
high-frequency data''}
\slink[doi]{10.1214/14-AOS1298SUPP} 
\sdatatype{.pdf}
\sfilename{aos1298\_supp.pdf}
\sdescription{This supplement contains technical proofs of the Lemmas
\ref{lem3}--\ref{lem4} as well as some interesting supplemental lemmas.}
\end{supplement}




\printaddresses
\end{document}